\documentclass[11pt, reqno]{amsart}
\usepackage{amsmath,amsthm,amssymb}

\theoremstyle{plain}
\newtheorem{thm}{Theorem}[section]
\newtheorem{prop}[thm]{Proposition}
\newtheorem{lem}[thm]{Lemma}
\newtheorem{cor}[thm]{Corollary}

\theoremstyle{definition}

\newtheorem{rem}[thm]{Remark}
\newtheorem{defn}[thm]{Definition}
\newtheorem{eg}[thm]{Example}
\newtheorem{subtitle}[thm]{}
\newtheorem{ex}{Exercise}[section]
\numberwithin{equation}{section}

\def\a{\alpha}

\def\d{\delta}
\def\D{\triangle}
\def\e{\epsilon}

\def\K{\nabla}
\def\l{\lambda}

\def\n{\,\vert\,}
\def\o{\theta}

\def\cd{{\mathcal{D}}}

\def\cg{{\mathcal{G}}}
\def\ch{{\mathcal{H}}}

\def\cm{{\mathcal{M}}}

\def\co{{\mathcal{O}}}

\def\li{\langle}
\def\ri{\rangle}
\def\n{\ \vert\ }
\def\tr{{\rm tr}}
\def\bs{\bigskip}
\def\ms{\medskip}

\def\ni{\noindent}
\def\ti{\tilde}
\def\p{\partial}

\def\Im{{\rm Im\/}}
\def\I{{\rm I\/}}

\def\diag{{\rm diag}}

\def\Gr{{\rm Gr}}
\def\sgn{{\rm sgn}}

\def\R{\mathbb{R} }
\def\C{\mathbb{C}}

\newcommand{\beg}{\begin{eg}}
\newcommand{\eeg}{\end{eg}}
\newcommand{\bthm}{\begin{thm}}
\newcommand{\ethm}{\end{thm}}
\newcommand{\bprop}{\begin{prop}}
\newcommand{\eprop}{\end{prop}}
\newcommand{\bcor}{\begin{cor}}
\newcommand{\ecor}{\end{cor}}
\newcommand{\blem}{\begin{lem}}
\newcommand{\elem}{\end{lem}}
\newcommand{\bca}{\begin{cases}}
\newcommand{\eca}{\end{cases}}
\newcommand{\brem}{\begin{rem}}
\newcommand{\erem}{\end{rem}}
\newcommand{\bpm}{\begin{pmatrix}}
\newcommand{\epm}{\end{pmatrix}}
\newcommand{\bbm}{\begin{bmatrix}}
\newcommand{\ebm}{\end{bmatrix}}
\newcommand{\bvm}{\begin{vmatrix}}
\newcommand{\evm}{\end{vmatrix}}
\newcommand{\bdefn}{\begin{defn}}
\newcommand{\edefn}{\end{defn}}
\newcommand{\bsub}{\begin{subtitle}}
\newcommand{\esub}{\end{subtitle}}
\newcommand{\bex}{\begin{ex}}
\newcommand{\eex}{\end{ex}}
\newcommand{\ben}{\begin{enumerate}}
\newcommand{\een}{\end{enumerate}}

\def\half{\frac{1}{2}}
\def\pt{\frac{\p}{\p t}}
\def\px{\frac{\p}{\p x}}
\def\py{\frac{\p}{\p y}}
\def\pz{\frac{\p}{\p z}}
\def\pzbar{\frac{\p}{\p \bar z}}
\def\pxi{\frac{\p}{\p \xi}}
\def\peta{\frac{\p}{\p \eta}}
\def\fh{\mathfrak{h}}

\begin{document}

\title[Space-time Monopole equation]
{On the space-time Monopole equation}
\author{Bo Dai}
\address{LMAM, School of Mathematical Sciences \\ Peking
University \\ Beijing 100871, P. R. China}
\email{daibo@math.pku.edu.cn}
\author{Chuu-Lian Terng$^*$}\thanks{$^*$Research supported in   part by NSF grant DMS- 0529756 }
\address{Department of Mathematics \\
University of California at Irvine, Irvine, CA 92697-3875}
\email{cterng@math.uci.edu}
\author{Karen Uhlenbeck$^\dag$}\thanks{$^\dag$Research supported in
part by Sid Richardson
Regents' Chair Funds, University of Texas system and NSF grant DMS-0305505\/}
\address{Department of Mathematics\\ University of Texas at Austin \\
Austin, TX 78712}
\email{uhlen@math.utexas.edu}

\ms


\begin{abstract}

The space-time monopole equation is obtained from a dimension reduction of the anti-self dual Yang-Mills equation on $\R^{2,2}$.  A family of Ward equations is obtained by gauge fixing from the monopole equation. In this paper, we give an introduction and a survey of the space-time monopole equation. Included are alternative explanations of results of Ward, Fokas-Ioannidou, Villarroel and Zakhorov-Mikhailov. The equations are formulated in terms of a number of equivalent Lax pairs; we make use  of the natural Lorentz action on the Lax pairs and frames.  A new Hamiltonian formulation for the Ward equations is introduced.  We outline
 both scattering and inverse scattering theory  and use B\"acklund transformations to construct a large class of monopoles which are global in time and have both continuous and discrete scattering data.  
 
\end{abstract}

\maketitle

\tableofcontents

\section{Introduction}\label{cq}

\ms

The self-dual Yang-Mills equations in $\R^4$ and their reduction
to  monopole equations in $\R^3$ have become central
topics of study and useful tools in modern geometry.  The same
self-dual equations in the case of a different signature of
$\R^{2,2}$ are not of the general type to be used much in
geometry. However, their dimensional reduction to the {\it space-time
monopole equation\/} in $\R^{2,1}$ yields an extremely interesting
system of non-linear wave equations which deserve to be better
known. These equations can be encoded in a Lax pair.  Moreover,
with a mild additional assumption and a gauge fixing they
can be rewritten for a map from $\R^{2,1}$ into a Lie group. These
equations differ only slightly from the usual wave map equation.

This article is meant to be an introduction to and a survey of the
literature on the space-time monopole equations. We also give a construction of the inverse scattering of the monopole equations via loop group factorizations.  These equations
form a hyperbolic system for a connection and a Higgs field, and
hence have a gauge symmetry.  A  simple restriction and coordinate
change produces the equation for a map into the gauge group.  This
last equation was introduced by Richard Ward, who studied them
using a version of Riemann-Hilbert problem  and twistor
theory. He produced the basic examples and a number of interesting
papers \cite{Wa88,Wa90,Wa95}. Hence the equation for the map is
referred to as either Ward's equation, or in his original
language, the modified chiral model. Additional work on the
equations is due to T. Ioannidou, W. Zakrewski
\cite{Io96,IoZa98a,IoZa98b}, Manakov and Zakharov \cite{MaZa81},
A. K. Fokas and Ioannidou \cite{FoIo01} and Villaroel \cite{Vi90}.
The last three references present both the continuous scattering
theory and the inverse scattering transform. The construction of a
complete set of soliton solutions has been carried out by the
first two authors in a previous paper \cite{DaTe04}.

The plan of the paper  is as follows. We derive the monopole
equations with their Lax pairs, paying special attention to the
difference between monopole equations in space and in space-time
in section 2. A family of Ward equations for maps into groups is
constructed in Section 3. Section 4 describes the action of the
Lorentz group on the Lax pair system and on frames. We make use of
Lorentz boosts in the construction of solitons and of the spacial
rotation group in deriving estimates in the appendix.  Next we
list special classes of solutions, so we can continue the
discussion with a lot of examples in mind. Section 6 contains a
very brief Hamiltonian formulation for the family of Ward
equations. In section 7, we introduce the transform which produces
the continuous scattering data as well as the inverse scattering
transform. Since the inverse scattering transform always exists,
this produces many global solutions to the equations that are decaying at spacial infinity. The details
of the fixed point theorem which yields continuous scattering data
for small initial data are in the appendix. In section 8 to 10, we  review B\"acklund transformations, and use these
transformations to construct soliton monopoles and monopoles with
both continuous and discrete scattering data.

Due to soliton theory, B\"acklund transformations and the inverse
scattering transform, we discover a very large class of solutions
which are global in time.  This is in contrast to the closely
related wave map equation from $\R^{2,1}$ to $G$.
It is a difficult theorem first of T. Tao \cite{Tao01}, extended by D.
Tataru \cite{Tat04}, to show that small initial data results in
solutions for all time. Whether the difference is entirely due to
integrability, or whether there is a deeper analytic theory  or
more examples to be found remains open.

\bs

\section{The Monopole equations}\label{cr}

\ms

The curvature $F_A= \sum_{i, j} F_{ij} \ dx_i \wedge dx_j$ of a
$u(n)$-valued connection $1$-form $A=\sum_{i=1}^4 A_i dx_i$ on
$\R^4$ is
$$
F_{ij}=[\K_i, \K_j] = -\p_{x_i} (A_j)+\p_{x_j} (A_i) + [A_i, A_j],
$$ where
$$
\K_i= \p_{x_i } - A_i.
$$ The connection $A$ is {\it anti-self-dual Yang-Mills\/}
(ASDYM) on $\R^4$  if
$$
\ast F_A= - F_A,
$$
where $\ast$ is the Hodge star operator with
respect to the Euclidean metric  $\sum_{i=1}^4 dx_i^2$ on $\R^4$.
The ASDYM on $\R^4$ written in coordinates
is
\begin{equation}\label{ai}
F_{12}= -F_{34}, \quad F_{13}= -F_{42},
\quad F_{14}= -F_{23}.
\end{equation}

The ASDYM has a Lax pair formulation.  The term
``Lax pair'' refers to any equation which is written as a ``zero
curvature" equation for a connection, or a portion of a
connection. This connection contains an additional complex
parameter $\mu$ is variously interpreted as a ``spectral'', ``twistor'', or
``Riemann-Hilbert parameter''.
Set
$$ z=x_1+ix_2, \quad w=
x_3+ix_4,
$$
$\K_z= \frac{1}{2}(\K_1-i\K_2)= \frac{\p}{\p z} - A_z$, $\K_{\bar
z}= \frac{1}{2} (\K_1+i\K_2)=\pzbar -A_{\bar z}$, and $\K_w,
\K_{\bar w}$ similarly. Since $A_i\in u(n)$,  $A_{\bar z} =
-A_z^*$ and $A_{\bar w}=-A_w^*$, where $B^*= \bar B^t$. The equation
\begin{equation}\label{ac}
[\K_{\bar w}+ \mu\K_z, \K_w -
\mu^{-1}\K_{\bar z}]=0
\end{equation}
is equivalent to the ASDYM \eqref{ai} on $\R^4$.  This is because \eqref{ac}
holds for all $\mu\in \C\setminus \{0\}$ if and only if the
coefficients of $\mu, 1$ and $\mu^{-1}$ of \eqref{ac} are zero,
which is \eqref{ai}.

If we assume the ASDYM connection $A$ is independent of $x_4$,
then $A_w=\half(A_t-i\Phi)$ and $A_{\bar w}=\half(A_t+i\Phi)$,
where $\Phi=A_4$ is the Higgs field, $A=A_t\ dt + A_z\ dz +
A_{\bar z}\ d\bar z$ is a connection $1$-form on $\R^3$.  Then $(A,\Phi)$
satisfies the Euclidean monopole equation (cf. \cite{AtHi88})
$$
D_A\Phi= \ast F_A,
$$
where $\ast$ is the Hodge star operator with respect to the metric
$dx^2+dy^2+ dt^2$. The Euclidean monopole equation is an important equation in  both geometry and physics (cf. 
\cite{Hi83,AtHi88}).

\ms The ASDYM on $\R^{2,2}$ is again $\ast F_A=- F_A$, but $\ast$
is the Hodge star operator with respect to the metric
$dx_1^2+dx_2^2 - dx_3^2 -dx_4^2$. In coordinates, the ASDYM is
\begin{equation}\label{aj}
F_{12}=F_{34}, \quad F_{13}=- F_{42}, \quad F_{14}= -F_{23},
\end{equation}
This equation has a Lax pair
similar to the ASDYM on $\R^4$,
\begin{equation}\label{aq}
[\K_{\bar w} +\mu \K_z, \ \K_w+\mu^{-1}\K_{\bar
z}]=0.
\end{equation}
(the only difference with \eqref{ac} is that the second operator
is taken with a plus sign).  In other words, $A$ satisfies
\eqref{aj} if and only if \eqref{aq} holds for all $\mu\in
\C\setminus \{0\}$.

We say $E(\cdots, \mu)$ is a {\it frame\/} of the Lax pair $[D_1(\mu), D_2(\mu)]=0$ if $D_1E=D_2E=0$.  

If $[D_1(\mu), D_2(\mu)]=0$ is a Lax pair of a PDE with spectral parameter $\mu$, then $[\ti D_1(\l), \ti D_2(\l)]=0$ is also a Lax pair of the same PDE with parameter $\mu= \frac{a\l +b}{c\l+d}$, where $\ti D_	i(\l)= \sum_{j=1}^2 f_{ij}(\mu) D_j(\mu)$ and $f_{ij}$ are meromorphic functions.  Moreover, if $E$ is a frame for the Lax pair $[D_1(\mu), D_2(\mu)]=0$, then  $F(\cdots, \l)= E(\cdots, \frac{a\l+b}{c\l+ d})$ is a frame for the Lax pair $[\ti D_1(\l), \ti D_2(\l)]=0$.  We say these Lax pairs are {\it equivalent\/}. 

If $D_1(\mu_0)g= D_2(\mu_0)g=0$,  then $g^{-1}[D_1(\mu), D_2(\mu)]g=0$ is also a Lax pair for the same PDE, which is called the Lax pair obtained from the Lax pair $[D_1(\mu), D_2(\mu)]=0$ by {\it fixing the gauge\/} at $\mu=\mu_0$.
 So Lax pairs of a PDE come in  many forms. We can choose an equivalent Lax pair and fix a gauge to make our computation easier or more transparent. For example, this is what we will do in later sections when we compute the $SO(2,1)$-actions on frames and  discuss the scattering and inverse scattering  of the monopole equations.   

By making a linear fractional transformations in $\mu$, we obtain an equivalent Lax pair. 
 Set $\mu= \frac{\l-i}{\l+i}$. Then
\begin{equation}\label{ba}
\bca2(\l+i)(\K_{\bar w}+\mu \K_z) =
L_1+iL_2,\\2(\l-i)(\K_w +\mu^{-1}\K_{\bar z})
=L_1-iL_2,\eca
\end{equation}
where
$$
\bca
L_1= \l(\K_{x_3} +\K_{x_1})
-(\K_{x_4}+\K_{x_2}), \\ L_2= \l(\K_{x_4}-\K_{x_2}) +
(\K_{x_3}-\K_{x_1}). \eca
$$
It follows that \eqref{aq} holds if and only if
$[L_1, L_2]=0$, i.e.,
\begin{equation}\label{ar}
[ \l(\K_{x_3}
+\K_{x_1}) -(\K_{x_4}+\K_{x_2}), \  \l(\K_{x_4}-\K_{x_2}) +
(\K_{x_3}-\K_{x_1})]=0
\end{equation}
holds for all $\l\in \C\setminus\{i, -i\}$. So  \eqref{ar} is an equivalent
Lax pair for the ASDYM on $\R^{2,2}$.  The use of $\mu$ in the section on scattering theory, and $\l$ in the rest of the literature on the Ward equation is unfortunately, confusing, but necessary. 

If we assume the connection $A$ is independent of $x_4$. then
$$
A_w =\frac{1}{2}( A_t - i\phi), \qquad  A_{\bar w} =
\frac{1}{2}(A_t + i\phi).
$$
where $A_{x_4} = \phi$ is now the Higgs field for our new {\it
space-time monopole equation\/}.  A calculation shows that we can
write the space-time monopole equation in the same form  as the
more familiar Euclidean monopole equation for a connection $A =
A_t\, dt + A_z\, dz + A_{\bar z}\, d\bar z$ and a Higgs field
$\phi$ (cf. \cite{Wa99}):
\begin{equation}\label{ad}  D_A \phi = * F_A.
\end{equation} Here we have used the $\ast$ operator for the
Lorentz metric $dx^2+dy^2-dt^2$.

In this paper, we always assume $(A, \phi)$ decays at spacial infinity. 

Equation \eqref{ad} has a  Lax pair induced from the Lax pair
\eqref{aq} for ASDYM:
\begin{equation}\label{ap}
\left[\half\K_t -\frac{i\phi}{2}
+\mu \K_z, \ \ \half\K_t +\frac{i\phi}{2} +\mu^{-1}\K_{\bar
z}\right]=0.
\end{equation}
Set
$$ \mu= \frac{\l-i}{\l+i}, \quad \xi=
\frac{t+x}{2}, \quad \eta=\frac{t-x}{2}.
$$
Then \eqref{ar} induces an equivalent  Lax pair for
the monopole equation in $(\xi, \eta, y)$
coordinates:
\begin{equation}\label{au}
\left[\l\K_\xi-\K_y +\phi, \  \l^{-1}\K_\eta-\K_y-\phi\right]=0.
\end{equation}
This is the Lax pair used by Ward.  We have:

\bprop The following statements are equivalent for a connection
$A$ on $\R^{2,1}$ and a Higgs field $\phi$:
 \ben
 \item $(A, \phi)$
is a solution of the space time monopole equation \eqref{ad} on
$\R^{2,1}$.
\item \eqref{ap} holds for all  $\mu\in \C\setminus\{
0\}$.
\item The linear system,
\begin{equation}\label{ao}
\bca(\frac{1}{2}\pt +\mu\pz)E= (\half(A_t+i\phi) + \mu A_z)E, \\
(\frac{1}{2}\pt + \mu^{-1}\pzbar)E = (\half(A_t-i\phi)
+\mu^{-1}A_{\bar z})E, \eca
\end{equation}
is compatible for complex parameter
$\mu$.
\item \eqref{au} holds for all $\l\in \C\setminus
\{0\}$.
\item The linear
system
\begin{equation}\label{ay}
\bca(\l\pxi -\py)F= (\l A_\xi -(A_y+\phi))F,\\
(\l^{-1}\peta-\py)F= (\l^{-1}A_\eta -(A_y-\phi))F. \eca
\end{equation}
is compatible for complex parameter $\l$.
\een
Moreover,  if $E(x,y,t, \mu)$ is a frame of \eqref{ap} (i.e., a solution of \eqref{ao}),
then
\ben
\item[(a)] $F(x,y,t,\l)= E\left(x,y,t,\frac{\l-i}{\l+i}\right)$ is a frame of \eqref{au} (i.e., 
a solution of \eqref{ay}), 
\item[(b)]  $E$ satisfies the reality condition
\begin{equation}
E(x,y,t,\bar\mu^{-1})^*E(x,y,t,\mu)=\I
\end{equation}
if and only if $F$ satisfies the reality condition
\begin{equation}
F(x,y,t,\bar\l)^* F(x,y,t,\l)=\I.
\end{equation}
\een
\eprop

\bs

\section{The Ward equations} \label{cs}

We call solutions  of linear
systems \eqref{ao} or \eqref{ay} that satisfy the reality condition   {\it monopole frames}.  Note that, unlike the case when the Lax pair is a full connection, locally there can be a serious lack of uniqueness in solving for a frame.  We resolve this lack of uniqueness away from $\mu\in S^1$ by observing that the spacial part of the Lax pair is a Cauchy-Riemann operator.  Frames, if they exist are unique if we require $E_\mu=\I$ at spacial infinity.  We expect the frames to exist at most points $\mu\not\in S^1$.  

When $\mu=\rho\in S^1$, the existence of frames is more problematic.  To obtain the Ward equation, we need extra assumptions, even for small initial data.

\bdefn Let $\rho\in S^1$, and  $(A,\phi)$ a solution  of the
space-time monopole equation such that $(A,\phi)$ decays at spacial infinity.  We say $(A,\phi)$ is {\it $\rho$-regular\/}, if there
is a smooth solution $k:\R^{2,1}\to U(n)$ such that
\begin{equation}\label{ak}
\bca(\frac{1}{2}\p_t+\rho \p_z)k= (A_{\bar w} + \rho A_z) k,\\
(\frac{1}{2}\p_t + \rho^{-1}\p_{\bar z})k= (A_w+\rho^{-1}A_{\bar
z})k,
\eca
\end{equation}
 and $k-\I$ and the first derivative of $k$ decays as $|z|\to
\infty$.  (Note that the  second equation of \eqref{ak} is the
Hermitian transpose of the first). \edefn

Let $f$ be a $U(n)$-valued map. Then$$f(\p_x-A)f^{-1} =
\p_x -(fAf^{-1} +(\p_x f)f^{-1})$$
is the gauge  transformation of $f$ on $\px -A$, or
$$f\ast A= fAf^{-1}+f_xf^{-1}.$$ Suppose $(A,\phi)$
is $\rho$-regular and $k$ is the solution of \eqref{ak}.  We fix
the gauge at $\mu=\rho$, i.e., we apply the gauge transformation
of $k^{-1}$ to the Lax pair \eqref{ap} to
get
\begin{equation}\label{am}
\left[\frac{1}{2}\p_t+\mu \p_z -
(\mu-\rho)\ti A_z, \ \half\p_t +\mu^{-1}\p_{\bar z}-
(\mu^{-1}-\rho^{-1})\ti A_{\bar z}\right]=0,
\end{equation}
where $\ti A_z= k^{-1}\ast A_z$ and $\ti A_{\bar z}= k^{-1}\ast A_{\bar z}$.
Or equivalently, the following linear system is compatible for an
open subset of parameters
$\mu$:
\begin{equation}\label{sb}
\bca(\half\pt +\mu\pz)E= (\mu-\rho) \ti A_z E,\\ (\half\pt
+\mu^{-1}\pzbar)E=(\mu^{-1}-\rho^{-1}) \ti A_{\bar z} E.
\eca
\end{equation}
Suppose $(A,\phi)$ is also $-\rho$-regular.  Then there exists
$g:\R^{2,1}\to U(n)$ satisfies \eqref{sb} with $\mu=-\rho$, i.e.,
\begin{equation}\label{al}
\bca (\half\p_t- \rho
\p_z) g= -2\rho \ti A_z g, \\ (\half\p_t -\rho^{-1}\p_{\bar z})g=
-2\rho^{-1}\ti A_{\bar z} g.
\eca
\end{equation}
A computation shows that
\begin{equation}\label{aw}
-(g_tg^{-1})_t + (g_xg^{-1})_x + (g_yg^{-1})_y + [g_tg^{-1}, \cos
\o\ g_xg^{-1}+ \sin \o\ g_yg^{-1}]=0,
\end{equation}
where $\rho=e^{i\o}$.  This is the one-parameter family of {\it
Ward equations\/} \cite{Wa88}.  We then obtain

\bprop
Suppose $E(x,y,t,\mu)$ is a frame for the solution of the
space-time monopole equation, (i.e., $E$ is a solution of the
linear system \eqref{ao}), and $E(x,y,t,\mu)$ are smooth at
$\mu=\pm e^{i\o}$.   Then $$g(x,y,t):=
E(x,y,t,e^{i\o})^{-1}E(x,y,t, -e^{i\o})$$ is a solution of the
Ward equation \eqref{aw}. \eprop

The Lax pair \eqref{au} is equivalent to the following Lax pair
\begin{equation}\label{bf} [\l\K_\xi -\K_y +\phi, \ \
\l(\K_y+\phi) -\K_\eta]=0.
\end{equation}
 We fix the gauge of \eqref{bf} at $\l=\infty$
(equivalent to fix the gauge of \eqref{ap} at $\mu=1$), i.e., take
the gauge transformation of $h^{-1}$ on \eqref{bf} to get
\begin{equation}\label{bc}
[\l \p_\xi -(\p_y -\hat A), \ \ \l
\p_y-(\peta- \hat B)]=0,
\end{equation}
where $\p_y h= (A_y -\phi)h$, $\p_\xi h =
A_\xi h$,  $\hat A=- 2h^{-1}\phi h$ and $\hat B= h^{-1}\ast A_\eta$. So
we have the following proposition:

\bprop The following statements are equivalent:
 \ben
 \item
Equation \eqref{bc} holds for all $\l\in \C$,
\item
$$\bca\p_\xi \hat
B= \p_y \hat A,\\ [\p_y-\hat A, \p_\eta -\hat B]=0. \eca
$$
\item  the
linear system
\begin{equation} \label{bd} \bca(\l\p_\xi-\p_y)H= -\hat
A H, \\ (\l\p_y - \p_\eta)H= -\hat B H,
\eca
\end{equation}
\een is locally solvable for an open subset of $\l\in \C$.  Moreover, if
$H(x,y,t,\l)$ is a solution of \eqref{bd} and is smooth at $\l=0$,
then $g=H(\cdots, 0)$ satisfies
\begin{equation}\label{be} \pt\left(\frac{\p g}{\p t}
g^{-1}\right) -\px\left(\frac{\p g}{\p x} g^{-1}\right)
-\py\left(\frac{\p g}{\p y} g^{-1}\right) -\left[\frac{\p g}{\p t}
g^{-1}, \frac{\p g}{\p x} g^{-1}\right]=0,
\end{equation}i.e., $g$
is a solution of Ward equation \eqref{aw} with
$\o=0$.
\eprop

 As a consequence of the above Proposition, we see that to construct solutions of the monopole equation that are $\l=\infty$-regular is suffices  to construct $H(\xi, \eta, y, \l)$ such that $(\l \p_\xi H-\p_y H)H^{-1}$ and 
 $(\l\p_y H- \p_\eta H)H^{-1}$ are independent of $\l$.

\bprop If a monopole is $\mu=\pm 1$ regular (i.e., $\l=\infty, 0$
regular), then it is gauge equivalent to a monopole $(A,\phi)$
such that $A_\xi=0$, $A_y=\phi$.  Conversely, if $(\hat A, \hat
B)$ satisfies \eqref{bc}, then $A_\xi=0$, $A_y=\phi= \hat A/2$,
and $A_\eta= \hat B$ is a monopole. \eprop

\bs

\section{The Action of $SO(2,1)$} \label{ct}

The Lorentz group $SO(2,1)$
is the group of all $g\in SL(3,\R)$ such that $g^t\I_{2,1}g=
\I_{2,1}$, where $\I_{2,1}=\diag(1,1,-1)$. The group $SO(2,1)$
acts on $\R^{2,1}$ by the standard action $g\cdot p= gp$ (here
$p\in\R^{2,1}$ is identified as a $3\times 1$ vector).  Given a
connection $A=A_1 dx+A_2 dy + A_3 dt$, a Higgs field $\phi$, and
$g\in SO(2,1)$, the action  $g\cdot (A,\phi)= (g\cdot A, g\cdot
\phi)$ is defined by $g\cdot A= g\cdot A_1\ dx + g\cdot A_2\ dy +
g\cdot A_3\ dt$, where
$$
(g\cdot A_i)(p)= A_i(g\cdot p), \quad
(g\cdot \phi)(p)= \phi(g\cdot p).
$$

The space-time monopole equation is invariant under the Lorentz
group $SO(2,1)$, i.e., if $(A, \phi)$ is a solution then so is
$g\cdot (A,\phi)= (g\cdot A, g\cdot \phi)$ for $g\in SO(2,1)$.
In order to make the scattering theory estimates tractable and to understand the $1$-solitons we need to understand the
 natural action of $SO(2,1)$ on frames
(solutions of linear system \eqref{ay}).   Since $SO(2)$ of the
$xy$-plane and $O(1,1)$ of the $xt$-plane generate $SO(2,1)$,  to
compute the explicit formula of the action of $SO(2,1)$ on frames,
it suffices to compute the action of the following one-parameter
subgroups on frames:
\begin{align*}
R(\o) &=\bpm \cos\o & -\sin\o&0\\ \sin \o &
\cos\o &0\\ 0&0&1\epm,\\T(s)&= \bpm \cosh s &0& \sinh s\\ 0&1&0\\
\sinh s& 0& \cosh s\epm,
\end{align*}
We also need the representation $\sigma: SO(2,1)\to SL(2,\R)$,
whose differential $d\sigma_e$ maps $e_{12}-e_{21}$ to $
-\half(e_{12}-e_{21})$, $e_{13}+e_{31}$ to $\half(e_{11}-e_{22})$,
and $ e_{23}+ e_{32}$ to $ -\half(e_{12}+ e_{21})$. So
$$\sigma(R(\o))= \bpm \cos\frac{\o}{2} &\sin\frac{\o}{2}\\ -\sin
\frac{\o}{2}& \cos \frac{\o}{2}\epm,\quad \sigma( T(s))= \bpm
e^{\frac{s}{2}}&0\\ 0& e^{-\frac{s}{2}}\epm.$$

The group $SL(2,\R)$ acts on the scattering parameter space $\C\cup\{\infty\}$ by the linear fractional
transformations:
$$
\bpm a& b\\ c& d\epm \ast \l= \frac{a\l +
b}{c\l+ d}.$$

\bthm \label{bu} The group $SO(2,1)$ acts on the Lax pairs and on the frames of the
space-time monopole equation.  If $F$ is a frame of $(A,\phi)$, then
$$(g\cdot F)(p, \l)=F(g\cdot p, \sigma(g)\ast \l)$$ is a frame for $g\cdot(A, \phi)$, where
$\sigma:SO(2,1)\to SL(2,\R)$ is the representation given above and
$\ast$ is the standard action of $SL(2,\R)$ on $\l$ via the linear
fractional transformation.
\ethm

\begin{proof} For the action of
$R(\o)$, it can be checked easily that if $E$ is a solution of
\eqref{ao} for $(A,\phi)$, then $\ti E$ solves \eqref{ao} for
$e^{i\o}\cdot (A, \Phi)$, where $$\ti E(z, t,\mu)=  E(e^{i\o} z,
t, e^{i\o}\mu).$$ (Here we identify $e^{i\o}$ with the rotation
matrix $R(\o)$). In other words,$$R(\o)\cdot E(p,\mu)= E(R(\o)p,
e^{i\o}\mu).
$$
To compute the action of $R(\o)$ on solutions $F$ of \eqref{ay},
we recall that $\mu=\frac{\l-i}{\l+i}$,
$\l=\frac{i(1+\mu)}{1-\mu}$, and$$F(p,\l)= E\left(p,
\frac{\l-i}{\l+i}\right),$$where  $E$ is a solution of \eqref{ao}.
But for $\ti \mu= e^{i\o}\mu$, the corresponding $\ti \l$ and $\l$
are related by$$\ti \l= \frac{i(1+e^{i\o})\l
-(1-e^{i\o})}{(1-e^{i\o}) \l + i (1+e^{i\o})} =
\frac{\cos\frac{\o}{2}\ \l + \sin \frac{\o}{2}}{-\sin \frac{\o}{2}
\ \l + \cos  \frac{\o}{2}} = \bpm \cos  \frac{\o}{2}& \sin
\frac{\o}{2}\\ -\sin  \frac{\o}{2}& \cos  \frac{\o}{2}\epm \ \ast
\l.
$$
In other words,
$$
(R(\o)\cdot F)(p, \l)= F(R(\o)p,
\sigma(R(\o))\ast \l).
$$

For the action of $T(s)$, we use the light cone coordinates $\xi,
\eta$. Then $T(s)$ maps $(\xi, \eta, y)$ to $(e^s\xi, e^{-s}\eta,
y)$.   Suppose $F(p,\l)$ is a solution of linear system \eqref{ay}
for $(A,\phi)$.  Given $s\in \R$, define$$\hat F(\xi, \eta, y,
\l)= F(T(s)\cdot p, e^s \l)= F(e^s \xi, e^{-s}\eta, y, e^s\l).$$
It can be checked easily that $\hat F$ is a frame for $T(s)\cdot
(A, \phi)$. But
$$
e^s\l= \bpm e^{\frac{s}{2}}&0\\ 0&
e^{-\frac{s}{2}}\epm\ \ast \l =\sigma(T(s))\ast \l,
$$
where $\ast$ is the standard action of $SL(2,\R)$ on $\C$.
\end{proof}

\bcor The group $SO(2,1)$ acts on frames $E$ of \eqref{ao} as
follows: If $E$ is a frame of $(A,\Phi)$, then $(g\cdot E)(p,
\mu)=E(g\cdot p, g\sharp \mu)$, where
\begin{align*} &R(\o)\sharp \mu =
e^{i\o}\mu,\\ &T(s)\sharp \mu = \frac{\cosh \frac{s}{2}\ \mu
+\sinh \frac{s}{2}}{\sinh \frac{s}{2} \ \mu +\cosh \frac{s}{2}}.
\end{align*}
\ecor

\bs

\section{Special classes of solutions}\label{cu}

Many examples of solutions to the monopole equations are obtained by
assuming additional conditions.  A first set  of examples come
from linear solutions of the wave equation.  If $H\subset G$ is any
abelian subgroup, the monopole equations for maps into $H$ are linear.
This provides us with a number of solutions to which we can later
apply B\"acklund transformations.
Let $\ch\subset \cg$ be the Lie subalgebra of $H$.
Suppose that $(A,\phi)= (a, a_0, \phi)\in \oplus^2 \fh\oplus \fh\oplus\fh$
decays at spacial infinity, where
$a_0=A_t, a=(a_1, a_2)= (A_x, A_y)$, and $\fh$ is the Lie algebra of $H$.
Then the monopole equations \eqref{ad} are written in space-time
(here $d$ and $\ast$ are spacial):
\begin{equation}\label{bg}
\bca\ast da= \frac{\p \phi}{\p t}, \\
\frac{\p a}{\p t} - da_0= \ast d\phi. \eca
\end{equation}
The gauge transformation of $e^{-u}$ appears as
$$
(a, a_0, \phi)\mapsto (a-du, a_0-\frac{\p u}{\p t}, \phi).
$$
A global way of fixing gauge suitably is to require that
$d\ast a=0$ (a spacial equation which is easily solved).
With this choice, we apply $d\ast$ to the second equation of
\eqref{bg} to obtain
$$
\pt(d\ast a) - d\ast d\ a_0 = d\ast\ast  d\phi = 0.
$$
Since $d\ast a=0$, $\D a_0=0$.  We assume $a_0$ decays at infinity,
hence $a_0=0$.  Finally, we conclude
\bprop The abelian monopole equations are equivalent to the linear
wave equation for $\phi$:
$$
\square \phi= (\p_t^2 - \p_x^2 -\p_y^2)\phi=0,
$$
with $\ast d\ a = \frac{\p \phi}{\p t}$, $d\ast a=0$, and $a_0=0$.
\eprop We make special note of the fact that the condition of
$\rho$-regular is not automatic even in the abelian case. It is
satisfied if we can find a gauge transformation $g=e^u\in H$ such
that for $\rho=e^{i\a}$ we have $D_1(\rho)g=D_2(\rho)g = 0$. But
$$
D_1(\rho)-D_2(\rho)= i(\sin \a \K_x -\cos\a \K_y -\phi),
$$
so we have $g^{-1}(D_1(\rho)-D_2(\rho))g = i(\sin\a \px -\cos\a \py)$,
which implies that
\begin{equation}\label{bh}
-\sin \a (a_1-\frac{\p u}{\p x}) +\cos \a (a_2 -\frac{\p u}{\p y})-\phi =0.
\end{equation}
Here $u$ decays at spacial infinity.
But \eqref{bh} is the ODE
$$
\frac{du}{d\tau}= a_1 \sin \a - a_2 \cos \a +\phi,
$$
where $d\tau= \sin \a d x -\cos \a d y$. The condition that $u$
decays at spacial infinity is given by integral conditions, so it
is not automatically true. However, if we can find a decaying solution of
\eqref{bh}, it is unique. Because all the variables satisfy wave
equations, if we find a solution $u$ for \eqref{bh} at time $t=0$,
we can propagate it using the wave equation to any $t$. It follows
that the condition of $\rho$-regular is a property of the initial
data alone.

Our second class of examples is well-known (cf. \cite{Uhl92}). It involves
solutions of the monopole equation which are invariant under time
translation.

\bprop
Suppose that $(A,\phi)$ is a solution of the monopole
equation, which is fixed under time translation. Then its Lax pair 
is gauge equivalent to the Lax pair of a harmonic map. \eprop

\begin{proof} The Lax pair \eqref{ap} becomes
$$
\left[\mu\K_z -\half(A_t+i\phi), \ \mu^{-1}\K_{\bar z} -\half(A_t-i\phi)\right]=0.
$$
Since the connection at $\mu=1$ is now a full connection in $\C$,
we fix the gauge at $\mu=1$, so the Lax pair is now a Lax pair for the harmonic map
$$
\left[\pz +\frac{1-\mu^{-1}}{2}g^{-1}(A_t+i\phi)g, \ \pzbar
+\frac{1-\mu}{2} g^{-1}(A_t-i\phi)g \right]=0,
$$
where $g:\C \to U(n)$ satisfies $$\bca \frac{\p g}{\p z} g^{-1}
=A_z +\half (A_t+i\phi), & \\ \frac{\p g}{\p \bar z} g^{-1}
=A_{\bar z} +\half (A_t-i\phi). & \\
\eca $$ 
Note that this differs from the Euclidean monopole
reduction, where the different sign produces Hitchin's self-dual
equation (\cite{Hi87}) rather than the harmonic map.
\end{proof}

\bcor  The Lorentz transformations of the stationary solution
corresponding to harmonic maps produce families of solutions to
the monopole equations.
\ecor

It can be checked that the solutions obtained in the above Corollary decay in space.
\ms

Special solutions of the monopole equations also come from the
reduction to $\R^{1,1}$, which involves the assumption that the
field $(A,\phi)$ are independent of one of the spacial variables,
say $y$. The solution of these equations will not directly yield
solutions of the monopole equations which decay in space, but we
will see that they do arise in the consideration of radially
symmetric solutions.

Assume $(A,\phi)$ is independent of $y$. Then the $A_y=\psi$
becomes a second Higgs field and the Lax pair \eqref{ap} reduces
to
$$
[\K_t -i\phi + \mu(\K_x+i\psi), \ \K_t +i\phi +\mu^{-1}(\K_x -i\psi)]=0.
$$The equations become
$$
\bca[\K_t, \K_x]=[\phi, \psi], \\
\K_t \psi =\K_x \phi,\\\K_t \phi= \K_x\psi.
\eca
$$
The usual Lax pair is obtained by restriction to characteristic
coordinates. Let $\K_t+\K_x= \K_\xi$, $\K_t-\K_x= \K_\eta$, and
$\phi_+=\phi+\psi$, $\phi_-= \psi-\phi$. We get three equations of
$\phi_+$ and $\phi_-$ encoded in the Lax pair
$$
[\K_\xi + \tau \phi_+, \ \K_\eta +\tau^{-1}\phi_-]=0.
$$
The wave map $\phi:\R^{1,1}\to G$ is obtained from this Lax pair
in the same manner that the Ward map is obtained from the Lax pair for
monopoles.

Finally we are interested in solutions with a radial symmetry. In
gauge theory, the symmetry is inspired by requiring that under the
pull-back of a space-time symmetry, the fields $(A,\phi)$ go to
gauge transformations of themselves.  For example, given $J\in
su(n)$ such that $e^{2\pi J}=\I$, then$$e^{i\o}\cdot (A,
\phi)(z,t)= (e^{J\o} A(e^{i\o}z, t) e^{-J\o}, \ e^{J\o}
\phi(e^{i\o}z, t) e^{-J\o})$$is an $SO(2)$-action. If the monopole
equation is invariant under this action, then the field
$$(A,\phi)=
(e^{J\o}ae^{-J\o}, e^{J\o}a_0e^{-J\o}, e^{J\o}\phi e^{-J\o}),$$
where $a, a_0$ and $\phi$ depend only on $r$ and $t$. Both $a$ and
$\phi$ must vanish at $r=0$ unless $J=0$.   We then make a
singular gauge transformation by $e^{J\o}$ to the form
$$
(A,\phi)=
(a_r, a_\o-J, a_0, \phi),
$$
where $J$ now indicates a  singularity
at $0$ and lack of suitable decay at $\infty$.

\bprop Fix a representation of $S^1\to U(n)$ given by
$e^{i\o}\mapsto e^{J\o}$. The monopole equations for a monopole
with radial symmetry induced by $J$ are equivalent to the
equations for $\K_r, \K_t, \phi$ and $\psi= \frac{a_\o-J}{r}$:
$$
\bca [\K_t, \K_r] = [\phi, \psi],\\ [\K_t, \psi] =[\K_r, \phi],\\
[\K_t, \phi]= [\K_r, \psi] +\frac{1}{r}\ \psi. \eca
$$
\eprop

\begin{proof} It is useful to notice that
\begin{align*}
\K_z &= \half(\K_x-i\K_y)= \frac{e^{-i\o}}{2}(\K_r -\frac{i}{r}
\K_\o),\\\K_{\bar z}&= \half(\K_x +i\K_y)=\frac{e^{i\o}}{2}(\K_r +
\frac{i}{r} \K_\o).
\end{align*}
 Substitute these expression in the
Lax pair \eqref{ap} and compare coefficients of $\mu, \mu^{-1}$
and the constant term to get$$[\K_t-i\phi, \ \K_r+
i\psi]=0,$$which yields two equations. The third equation come
from the equation
$$
[\K_t-i\phi, \K_t+i\phi]
+[e^{-i\o}(\K_r-\frac{i}{r}\K_\o), \ e^{i\o}(\K_r +\frac{i}{r}
\K_\o)]=0,
$$
This last equation has an extra term which spoils the
integrability of the system, and prevents the three equations from
being encoded as a Lax pair.
\end{proof}

\bs

\section{Hamiltonian structures} \label{cv}

Let $(p, q)$ denote the
standard variables for the cotangent bundle $\cm$ of the space of
rapidly decaying maps from $\R^2$ to $U(n)$, and $H:\cm\to \R$ the
functional defined by
$$
H=\half \int\int_{\R^2} ||p||^2 + ||q^{-1}dq||^2\ dx dy.
$$  We introduce a Hamiltonian formulation of the Ward equation, which does not seem to have appeared in the literature before.   For each unit direction $v$ of $\R^2$, (i)  the
symplectic form $w^v$ is the sum of the standard symplectic
form on cotangent bundle $\cm$ and an extra closed $2$-form
depending on $v$, (ii) the Hamiltonian system of $H$
with respect to $w^v$ is the Ward equation with $\rho= v$.     

Let $\Theta$ be the canonical $1$-form on $\cm$ defined by
$$
\Theta_{(p,q)}(\d p, q^{-1}\d q) = p(q^{-1}\d q).
$$
The standard symplectic form on $\cm$ is $d \Theta $. We identify
the tangent and cotangent spaces in the Lie algebra formulation
with the $L^2$ inner product
$$
\li A, B\ri =\int_{\R^2} \tr(AB)\ dx
dy.
$$
Use the  Cartan formula
$$
d \Theta(\xi_1, \xi_2)= \xi_1 (\Theta(\xi_2))-\xi_2(\Theta(\xi_1))
-\Theta([\xi_1, \xi_2])
$$
to compute $d \Theta $  to get
\begin{align*} &d \Theta_{(p,q)}((\d_1p,
q^{-1}\d_1 q), (\d_2 p, q^{-1}\d_2 q))\\ &\quad = \li \d_1 p,
q^{-1}\d_2 q\ri -\li \d_2 p, q^{-1}\d_1 q \ri -\li p, [q^{-1}\d_1
q, q^{-1} \d_2 q]\ri.
\end{align*} A computation shows that the
following $2$-form is closed:
$$
\tau^v_{(p,q)}((\d_1p, q^{-1}\d_1
q), (\d_2 p, q^{-1}\d_2 q)) = \li q^{-1}q_v, [q^{-1}\d_1 q,
q^{-1}\d_2 q]\ri,
$$where $q_v= dq(v)= q_x \cos \o + q_y \sin \o$
if $v=e^{i\o}$. Then  $w^v= d\Theta + \tau^v$ is a symplectic form
on $\cm$.  In fact,
\begin{align*}
&w^v_{(p,q)}((\d_1p, q^{-1}\d_1 q), (\d_2 p, q^{-1}\d_2
q))\\&\quad = \li \d_1 p, q^{-1}\d_2 q\ri -\li \d_2 p, q^{-1}\d_1
q \ri +\li -p + q^{-1}q_v, [q^{-1}\d_1 q, q^{-1} \d_2 q]\ri
\end{align*}
 Since
 $$
 dH_{(p,q)} (\d p, \d q)= \li \d p,
p\ri -\li q^{-1} \d q, d^*(q^{-1}dq)\ri,$$a direct computation
shows that the Hamiltonian flow for $H$ with respect to the
symplectic structure $w^v$ is
$$
\bca q^{-1}q_t= p,\\ p_t +[-p+q^{-1}q_v, q^{-1}q_t]= d^*
(q^{-1}dq). \eca
$$
So we have proved

\bprop  Given a unit vector $v=e^{i\o}$ in $\R^2$, the $2$-form
$w^v$ is a symplectic form, and  the Hamiltonian system of $H$
with respect to $w^v$ is the Ward equation \eqref{aw}.
\eprop

In \cite{IoZa98b}, Ioannidou and Zakrzewski considered another
family of Ward equations, and studied the Lagrangian and
Hamiltonian formulations for that family of equations.

\bs

\section{Scattering theory\/}\label{cw}

Scattering theory for the Ward equations has been treated by a
number of authors, including Manakov and Zakharov \cite{MaZa81},
Villarroel \cite{Vi90}, and Fokas and Ioannidou \cite{FoIo01}. We
include a brief of synopsis and interpretation of the results.  In particular, we construct the inverse scattering transform via loop group factorizations.

Ward's original analysis of the space-time monopole equations is
via twistor theory. We recognize features of this analysis in what
follows, although we do not go into the twistor formulation.
Recall that
$$
D_1(\mu)= \half \K_t
-\frac{i\phi}{2} +\mu \K_z, \quad D_2(\mu)= \half \K_t
+\frac{i\phi}{2} +\mu^{-1}\K_{\bar z}.
$$
We rewrite the Lax pair
as the linear system consisting of first a spacial operator
\begin{align}\label{ca}
D_s(\mu) & = D_1(\mu)-D_2(\mu) = \mu \K_z -\mu^{-1}\K_{\bar z} - i\phi\\
&=\frac{(\mu-\mu^{-1})}{2}\K_x -i\frac{(\mu+\mu^{-1})}{2}\K_y -
i\phi.
\end{align}
The second operator in the Lax pair we write
as one of a family of time operators
\begin{align*}D_t(\mu)&=
D_1(\mu)+D_2(\mu) +\a D_s(\mu)\\
&= \K_t +\half(\mu+\mu^{-1})\K_x -
\frac{i}{2}(\mu-\mu^{-1})\K_y + \a D_s(\mu).
\end{align*} The
operators can be rescaled (i.e., multiply by a scalar function of
$\mu$), so the points $\mu=0$ and $\mu=\infty$ are included by
changing the scaling factor (for example, if we multiply by $\mu$ to
the operators, then they are defined at $\mu=0$, and if we
multiply by $\mu^{-1}$ to the operators, then they are defined at
$\mu=\infty$).  Note that $D_1E=D_2E=0$ if and only if $D_sE=
D_tE=0$.

If $\mu\not\in S^1$, then the spacial part of the connection,
$D_s(\mu)$, can be thought of as containing a $\bar\p$ operator in
the complex structure on $\R^2$ given by the complex coordinate
$$w=\frac{\mu^{-1} z+\mu\bar z}{a(\mu)}, \quad \bar w=
\frac{\bar\mu^{-1}\bar z +\bar \mu z}{\bar a(\mu)}.
$$
(In fact, $D_s(\mu)= \K_{\bar w}$). The factor $a(\mu)$ does not
change the complex structure. Note at $\mu=\infty$, $w=\bar z$,
and at $\mu=0$, $w= z$.

Suppose the Higgs field $\phi$ and the connection $A$ decay at
infinity. For fixed time and every $\mu\in
\C\cup\{\infty\}\setminus S^1$ the connection $D_s(\mu)$
determines a bundle holomorphic in the complex parameter
$w=w(\mu)$ on $\R^2\cup \{\infty\}= S^2$. By a Theorem of
Grothendieck, this bundle is the sum of line bundles
$$
L_\mu=
L_1\oplus L_2\oplus \cdots \oplus L_n
$$
with first Chern classes $c(1)\leq c(2)\leq \cdots \leq c(n)$. The
reality condition insures that the Chern classes at $\bar\mu^{-1}$
are the negative of the Chern classes at $\mu$, since the bundle
$L_{\bar \mu^{-1}}$ is dual to the bundle $L_\mu$ at $\mu$.  We
call this sequence $\vec{c}(\mu)= (c(1), \ldots, c(n))$ the {\it
Chern vector at $\mu$}.

\bthm
Suppose we have a solution of the monopole equation in a time
interval $[T_1,T_2]$. Then for every $\mu\in
\C\cup\{\infty\}\setminus S^1$, the Chern vector of the spacial
holomorphic bundle is preserved under the flow  in
time.
\ethm

\begin{proof} As a warm-up, we first prove this for
$\mu=0$.  Then the spacial connection is $\K_{\bar z}$ and we
choose $D_1(0)= \K_t-i\phi$ as the evolution operator.  We have
$[\K_t-i\phi, \ \K_{\bar z}]=0$.  We may make a complex gauge
transformation so that $D_1(0)= \K_t -i\phi= h^{-1}\circ \pt \circ
h$. Hence the gauge equivalent Lax pair gives
$$
\pt (h\circ
\K_{\bar z}\circ h^{-1})=0.
$$
The complex structure on the bundle determined by $\K_{\bar z}$ is
carried into a structure which is a gauge equivalent one. Hence
they have the same splitting, and so the same Chern vectors.

A similar computation occurs for each $\mu$, where we make the
choice of the time direction at $\mu$ to be real.  For simplicity,
let $\mu= e^{i\o}\rho$, where $\rho\in \R^+$ and rotate variable
by $z\mapsto e^{i\o}z$ so that in this coordinate system
$\mu=\rho$ is real.  Choose as a suitable evolution operator
$$
D_1(\rho)+ \rho^2 D_2(\rho) =
\frac{1+\rho^2}{2}\K_t +\rho \K_x -\frac{i(1-\rho^2)\phi}{2}.
$$
The derivatives which appear are in the direction of
$\frac{\p}{\p\tau} =\frac{1+\rho^2}{2} \pt + \rho\px$. We can
again make a gauge transformation so that
\begin{align*} &\K_\tau -
\frac{1-\rho^2}{2}\ i\phi = h^{-1}\circ\frac{\p}{\p \tau}\circ
h,\\& \frac{\p}{\p \tau} (h\circ D_s(\rho)\circ h^{-1})=0.
\end{align*}
Hence the complex structure of the bundle does not changes under a
flow in the $\tau$ variable. However, time translation is a
translation in the $\tau$ direction followed by a translation in
the $x$ variable.  Translation in $S^2$ by $x$ is holomorphic and
does not change the complex structure. Hence the splitting type of
$L_\mu=L_\rho$ does not change under the flow in time.  We
conclude that the Chern vector $\vec{c}(\mu)$ is preserved.
\end{proof}

What do we expect? Since the space of $\bar \p$ derivatives which
lead to non-trivial splitting has codimension at least two, we
expect that for most choices of initial data, the non-trivial
splitting occur at isolated points.  At these points in $\C\setminus S^1$, we expect to have singularities (poles and zeros).  We expect continuous
scattering data to be defined as a jump across $\mu\in S^1$.

\bdefn A rapidly decaying spacial pair $(A,\phi)$ is said to have
{\it continuous scattering data\/} if the frame, which
solves
$$
\bca D_s(\mu)E_\mu=\left(\frac{\mu-\mu^{-1}}{2}\K_x -
\frac{i(\mu+\mu^{-1})}{2}\K_y -i\phi\right)E_\mu =0,\\
E_\mu(\infty)=\I, \quad E_{{\bar u}^{-1}}=
(E_\mu^*)^{-1},\eca
$$
has solutions $E_\mu^\pm$, which are holomorphic in a 
$$\co^{\pm}_\e=\{\mu\in \C\n 1<|\mu|^{\pm 1}<1+\e\}$$ for some $\e>0$. Moreover, we assume that the limits
$$
\lim_{\mu\in \co^\pm, \mu\to e^{i\o}} E_\mu =
S_\o^\pm
$$
exist. It follows from the reality condition that
$S_\o^-= {(S_\o^+)^*}^{-1}$. We call the non-negative Hermitian
matrix
$$S_\o=(S_\o^{-})^{-1}S_\o^+ = (S_\o^+)^\ast S_\o^+$$
the
{\it scattering matrix\/}.
\edefn

Let $W^{2,1}$ denote the  space of maps $f$ whose partial derivatives up to second order are in $L^1$.

\bprop Assume there is a gauge in which $(A,\phi)$ is rapidly
decaying in spacial variables, moreover, assume $(A,\phi)$ is
small in $W^{2,1}$. Then the Chern vector $\vec{c}(\mu)=0$ at
every $\mu\in \C\cup\{\infty\}\setminus S^1$. Moreover, the
continuous scattering matrix $S_\o$ exists, $\I-S_\o$ decays for
each $\o$, and the scattering matrix $S_\o$ satisfies
\begin{itemize}
\item[(a)] $\I-S_\o$ is small in
$L^\infty$,
\item[(b)] $S^\ast_\o = S_\o \geq 0$,\item[(c)]
$d_{s,\o}S_\o= (-\sin\o \px +\cos\o
\py)S_\o=0$.
\end{itemize}
\eprop

\begin{proof} The existence of
$E_\mu$ away from the circle $\mu\in S^1\subset \C P^1$ is a
straight forward iteration argument involving estimates in
$L^\infty$ on $E_\mu$ using $(A,\phi)$ small in $L^1\cap
L^\infty$.  We relegate the estimate as $\mu\to e^{i\o}\in S^1$ to
an appendix.  Once the basic estimate are in place, proof of
regularity and holomorphic dependence on $\mu$ are straight
forward.  We explicitly derive (a) in the appendix.  To obtain
(b), note that $(E_{\bar u^{-1}})^{-1}= E_\mu^\ast$ due to the
reality condition on $(A,\phi)$.  Hence $(S_\o^-)^{-1}=
(S_\o^+)^\ast$.

To obtain (c), notice that
$$
D_s(e^{i\o})= -i (d_{s,\o}
-\psi),
$$
where $\psi= -\psi^* = \sin\o A_x -\cos\o A_y +\phi$.
Since $D_s(\mu)E_\mu=0$, it follows from the definition of $S^+$
that $D_s(e^{i\o})S^+=0$, so $d_{s,\o}S^+=\psi S^+$. But
$$
d_{s,\o}((S^+)^*)= (d_{s,\o}S^+)^\ast = (\psi S^+)^\ast=
(S^+)^* \psi^*=-(S^+)^* \psi.
$$
Now
compute
\begin{align*}d_{s,\o}S&= d_{s,\o}((S^+)^*S^+) =
(d_{s,\o}(S^+)^*)S^+ + (S^+)^*d_{s,\o}S^+\\& = -(S^+)^*\psi S^+ +
(S^+)^*\psi S^+=0,
\end{align*}
and (c) follows.
\end{proof}

\ms

\bcor Assume $(A,\phi)$ is a smooth solution in $\R^2\times (T_1,
T_2)$ and decays in spacial variables, and has a smooth continuous scattering
data. Then
$$
0= \left(\pt +\cos\o \px +\sin \o \py
\right)S_\o.
$$
\ecor

The Corollary assumes that the scattering theory is differentiable
in $t$.  Since $D_t(\mu)=D_1(\mu)+D_2(\mu)= \K_t +\mu\K_z
+\mu^{-1}\K_{\bar z}$ has the property that $[D_s(\mu),
D_t(\mu)]=0$, we conclude that $D_t(\mu)E_\mu=0$.  Now the result
follows by the same method as (c) in the Proposition.

\bcor
If $(A,\phi)$ and $(\ti A, \ti\phi)$ are gauge equivalent,
then they have the same scattering data.
\ecor

\begin{proof} Suppose
$(\ti A, \ti\phi)$ is the gauge transformation of $(A,\phi)$ by a
unitary map $u$. If $E_\mu(x,y,t)$ is the frame for $(A,\phi)$,
then $u(x,y,t)E_\mu(x,y,t)$ is the frame for $(\ti A, \ti \phi)$.
Hence the limits $\ti S^+_\o = uS^+_\o$ and $\ti S^-_\o= u S_\o^-
$. But $u$ is unitary, so $\ti S=S$.
\end{proof}
Inverse scattering theory is simpler than scattering theory.  To
dispose the gauge ambiguity, we need to make a choice somewhere.
We choose to do this at $\mu=1$.  First we note that if we have
initial scattering data $S_\o(x,y)$ satisfying
$d_{s,\o}S_\o(x,y)=0$, then the scattering data $S_\o(x,y,t)$ for
the solution at time $t$ should satisfy
$$\frac{\p S_\o}{\p t}
+\cos\o \frac{\p S_\o}{\p x}+\sin\o \frac{\p S_\o}{\p y}=0. $$
Since $d_{s,\o}S_\o(x,y)=0$, there exists $s$ such
that
$$S_\o(x,y)= s(x\cos\o + y\sin\o ,\o).$$
The linear evolution
equation for $S_\o(x,y,t)$ implies that
$$
S_\o(x,y,t)= s(x\cos\o +y \sin\o -t,\o).
$$
The problem is now to write
$$S_\mu(x,y,t)=
E_\mu^-(x,y,t))^{-1} E^+_\mu(x,y,t), \quad \mu\in S^1,
$$ where
$E_\mu^-$ extends holomorphically to $\mu$ inside the unit circle,
$E_\mu^+$ extends holomorphically outside the unit circle, and
$E^-_{\mu}= ((E^+_{\bar\mu^{-1}})^*)^{-1}$. We can always do this,
and the ambiguity corresponds to the gauge transformations.  The
condition that $S_1(x,y,t)=\I$ is equivalent to the solution being
$1$-regular.  However, we choose a method of factoring which yield
a unique solution for all scattering data. The inverse scattering
was given in \cite{Vi90} by Villarroel and in \cite{FoIo01} by
Fokas-Ioannidou. We prove it using the Iwasawa loop group
factorization theorem of Pressley and Segal \cite{PrSe86}. Since
$S_\o$ is Hermitian symmetric and non-negative, there is a
Hermitian symmetric matrix $P_\mu$ such that $S_\o= P_\mu^2$,
where $\mu= e^{i\o}$.   The Iwasawa loop group factorization of
Pressley and Segal \cite{PrSe86} says that we can factor $$P_\mu=
U_\mu E_\mu^+$$ uniquely such that $ U_\mu$ is in $U(n)$, $U_1=
\I$, and $E_\mu^+$ extends holomorphically to outside the unit
circle.  Set $E^-_\mu= ((E^+_{\bar\mu^{-1}})^*)^{-1}$. Now $S_\o=
(E_\mu^-)^{-1}E_\mu^+$ when $\mu=e^{i\o}$. Since $S_\mu$ is
smooth, $E_\mu^\pm$ is smooth.

\bthm \label{bz} (Villarroel \cite{Vi90}, Fokas-Ioannidou
\cite{FoIo01})
\par Let
$S_\o(x,y,t)= s(x\cos\o + y\sin\o -t,\o)$.  Then the factorization
described above
$$
S_\o= (E_\mu^-)^{-1}(x,y,t)
E^+_\mu(x,y,t),
$$
where $\mu=e^{i\o}$ and $E_\mu^\pm$ can extend
holomorphically to $|\mu|^{\pm 1} >1$ yields smooth frame for a
solution of the space-time monopole equation. \ethm

\begin{proof} We
need to show that $E_\mu^\pm(x,y,t)$ generates a solution to the
monopole equation. To do this, note by construction that
$E_\mu^\pm$ is holomorphic in $|\mu|^{\pm 1} > 1$.  The operator
$d_{s,\o}$ is a directional derivative, $d_{s,\o}S=0$ and
$S=(S^-)^{-1}S^+$, so
$$
0=d_{s,\o}S= (d_{s,\o}((S^-)^{-1}))S^+ +
(S^-)^{-1}d_{s,\o}S^+.
$$
Thus we have
$$ (d_{s,\o}S^+)(S^+)^{-1} =-S^- d_{s,\o}((S^-)^{-1})
=(d_{s,\o}S^-)(S^-)^{-1}. $$
Note that
$$
d_{s,\mu}=
\frac{i(\mu-\mu^{-1})}{2} \ \px + \frac{\mu+\mu^{-1}}{2}\ \py
$$
is the meromorphic extension of $d_{s,\o}= d_{s,e^{i\o}}$. So by
meromorphic extension, we obtain the identity on $0<|\mu| <\infty$
$$
(d_{s,\mu}E_\mu^+) (E_\mu^+)^{-1} = (d_{s,\mu} E_\mu^-)
E_\mu^{-1}.
$$
Since the left hand side is meromorphic in $|\mu|>1$ with a simple
pole at $\infty$ and the right hand side is meromorphic on
$|\mu|<1$ with  a simple pole at $0$, both sides are meromorphic
in $\C$ with simple poles at $0$ and $\infty$.  Let $E_\mu
=E_\mu^\pm$ when $|\mu |^\pm >1$. Then $E_{\bar \mu^{-1}}=
(E_\mu^*)^{-1}$. It follows that
$$
\xi_{\mu}:= (d_{s,\mu}E_\mu)
E_\mu^{-1}= \mu C_1 +\mu^{-1} C_{-1} + C_0
$$
for some $C_i(x,y,t)$. But $d_{s,\bar\mu^{-1}}^*= d_{s,\mu}$ and
$E_{\bar\mu^{-1}}= (E_\mu^*)^{-1}$, so $\xi_{\bar\mu^{-1}}^*=
-\xi_\mu$. Hence $C_{-1}= -C_1^*$ and $C_0^*= -C_0$. Write $C_1=
A_z$, $C_{-1}= A_{\bar z}$ and $C_0=\phi$, then we have
 $$
(d_{s,\mu} E_\mu)E_\mu^{-1}= \mu A_z - \mu^{-1}A_{\bar z} +\phi,
$$
where $A_z=-(A_{\bar z})^*$ and $\phi^*=-\phi$.  In other
words, we have proved $(D_1(\mu)-D_2(\mu)) E_\mu=0$.

The proof of the evolution equation, i.e., $D_t(\mu) E_\mu=0$, is
similar and we do not carry out here.
\end{proof}

Of course, if we start with
an initial condition which has only continuous scattering data, we
will not necessarily obtain the same initial data by the inverse
scattering transform, but obtain a gauge equivalent
solution.

\bcor Let $(A,\phi)$ be a solution of the space-time monopole
equation rapidly decaying in the spacial variables. Assume in
addition that the Chern vector $\vec{c}(\mu)=0$ for all $\mu\in\C
\setminus S^1$. Assume also that the solution has continuous
scattering data for all $t$. Then the solution obtained from the
scattering data at $t=0$ agrees with the given solution up to
gauge transformation.
\ecor

\begin{proof} We know from the assumptions that
$S_\o(x,y,t)= s(x\cos\o +y\sin \o -t,\o)$, since $S_\o$ is unique,
$$S_\o= \lim_{\mu\to e^{i\o}, |\mu|>1} (E_{\bar\mu^{-1}})^*
E_\mu.
$$
Hence the frame provides a factorization. This
factorization is unique up to a unitary matrix $u=u(x,y,t)$. This
unitary matrix gives a gauge transformation between the original
solution and the one constructed by inverse
scattering.
\end{proof}

\bs

\section{$1$-soliton monopoles} \label{cx}

In addition to continuous scattering data, solutions of monopole equation may also have discrete scattering data.  We first construct monopoles whose frames have one simple pole, in later sections we construct frames with multiple poles and show how to combine them with continuous scattering data.  

The building blocks of the discrete scattering data are one-solitons, which are easy to describe.  We have the harmonic maps $\phi:\R^2\cup\{\infty\}\to SU(n)$, which yield time-independent solutions to the Ward equation.  Among these, we have one-unitons, which come from holomorphic maps into Grassmannians. We also have the Lorentz transformations of these stationary one unitons.  This family makes up the one-solitions. It is somewhat more difficult to show that every one-soliton, defined in terms of a single pole for the frame, is of this type.

We need to use another gauge equivalent Lax pair to construct soliton solutions.  
 If a monopole $(A,\phi)$ is $\l=0$
regular, then  we can fix the gauge of \eqref{au} at $\l=0$ to get
$$
[\l(\p_\xi -\ti A_\xi)- \p_y,\  \l(\p_y +2\ti\phi)-\p_\eta]=0,$$
where $\ti\phi= -\ti A_y$, and $\ti A_\eta=0$.  The above Lax pair
is equivalent to
\begin{equation}\label{cf} [\tau\p_y
-\p_\xi +\ti A_\xi, \ \tau\p_\eta -\p_y - 2\ti\phi]=0.
\end{equation}
(Here $\tau= \l^{-1}$).

The spectral parameters $\mu, \l, \tau$ in Lax pairs \eqref{ap},
\eqref{au}, and \eqref{cf} are related by
$$\mu= \frac{\l-i}{\l+i}, \quad \tau=\l^{-1},\quad {\rm so\ \
}\tau=\l^{-1} = \frac{i(\mu-1)}{\mu+1}.
$$

The above discussion gives the following Proposition (cf.
\cite{Wa88}):

\bprop \label{bo} Suppose there is a smooth $GL(n,\C)$-valued $\psi(x,y,t,\tau)$
defined for $(x,y,t)\in \R^{2,1}$ and $\tau$ in an open subset of
$\C$ such that
\ben
 \item $P:=(\tau\p_y \psi- \p_\xi \psi)\psi^{-1}$
and $Q:=(\tau \p_\eta \psi- \p_y \psi)\psi^{-1}$ are independent
of $\tau$,
 \item $\psi(x,y,t,\bar\tau)^*\psi(x,y,t,\tau)=\I$.
 \een
Let $A$ be a connection and $\phi$ a Higgs field defined by
$A_\xi=- P$, $A_\eta=0$, and $A_y= -\phi= - Q/2$.   Then
$(A,\phi)$ is a solution of the monopole equation.   Conversely,
every solution of the monopole equation that is regular at $\l= 0$
is gauge equivalent to a solution of this type. \eprop

\bdefn  \label{co}
A map $\psi$ satisfies (1) and (2) of Proposition \ref{bo} is called a {\it Ward frame\/} if
$P, Q$ decay in spacial infinity and $\psi(x,y,t,\infty)=\I$.  A Ward frame is a {\it Ward soliton frame\/} if $\psi$ is rational in $\tau$.
 \edefn
 
\bdefn  A solution $(A,\phi)$ of the monopole equation is called a {\it $k$-soliton\/} if  it is regular at $\mu=-1$ and has a monopole frame $E$ that is rational in $\mu$ with $k$ poles counted with multiplicity. \edefn

Note that a Ward soliton frame is a monopole frame.

If $f:S^2\to GL(n,\C)$ is rational with one simple pole at $\tau=\a$, then it can be checked that (cf. \cite{Uh89})  $f$ must be of the form 
$$g_{\a, \pi}(\tau)= \I + \frac{\a-\bar\a}{\tau-\a}\ \pi^\perp,$$ where $\pi^\perp= \I-\pi$ and $\pi$ is a Hermitian projection $\pi$ of $\C^n$.  We identify the space of rank $m$ Hermitian projections of $\C^n$ as $\Gr(m,\C^n)$ via the map $\pi\mapsto \Im(\pi)$.   So a  Ward $1$-soliton frame must be of the form
$g_{\a,\pi(x,y,t)}(\tau)=\I+ \frac{\a-\bar\a}{\tau-\a}\ \pi^\perp(x,y,t)
$ for some constant $\a\in
\C\setminus \R$ and $\pi:S^2\times \R\to \Gr(k, \C^n)$.  Note that
$(\tau\p_y \psi-\p_\xi\psi)\psi^{-1}$ and $(\tau\p_\eta\psi-
\p_y\psi)\psi^{-1}$ are independent of $\tau$ if and only if the
residue at $\tau=\a$ is zero.  This implies

\bprop \label{bw} Given $\a\in \C\setminus \R$ a constant and
$\pi:\R^{2,1}\to \Gr(k,\C^n)$ a smooth map, then $g_{\a,\pi}(\tau)=\I+ \frac{\a-\bar\a}{\tau-\a}\ \pi^\perp$ is a  Ward soliton frame if and only
if
\begin{equation}\label{so}
\bca (\a\p_y\pi -\p_\xi \pi)\pi=0, & \\
(\a\p_\eta\pi -\p_y\pi)\pi=0. & \eca
\end{equation}
Moreover, if $\pi$
is a solution of \eqref{so}, then there exists a holomorphic map
$\pi_0:S^2\to \Gr(k,\C^n)$ such that $\pi(x,y,t)=
\pi_0(y+\a\xi+\a^{-1}\eta)$. \eprop  

Note that if $\a=\pm i$, then $x+\a\xi +\a^{-1}\eta= y\pm i x$ and
$\pi(x,y,t)=\pi_0(y\pm ix)$.  So the $1$-soliton $g_{\pm i, \pi}$ is  a
$1$-uniton harmonic map. 

The  $SO(2,1)$-actions  described in section 4 of $1$-unitons are  $1$-soliton monopoles.  In fact, we have

\bprop \label{cm}
Given $\a=re^{i\o}\in \C\setminus\R$, let $e^s= r$,
$e^c=\csc\o + \cot\o= \cot(\o/2)$, and $h= T(c)R(-\pi/2) T(s)\in
SO(2,1)$, where $T(c)$ and $R(\o)$ are $1$-parameter subgroups of
$SO(2,1)$ defined in section 4. Let $\pi_0:S^2\to \Gr(k, \C^n)$ be
a holomorphic map. Then the action of $h$ on $1$-uniton frame
$g_{i,\pi_0}$ gives rise to a monopole solution that is gauge
equivalent to the $1$-soliton given by $g_{\a,\pi}$, where $\pi(x,y,t)= \pi_0(\frac{i}{\sin\o} (y+\a\xi +\a^{-1} \eta))$.  In
other words, all $1$-solitons monopoles are obtained from the action of
$SO(2,1)$ on $1$-unitons up to gauge equivalence. \eprop

\begin{proof} Recall that $\l=\tau^{-1}$, and $$F(x,y,t,\l)=
g_{i,\pi_0(y+ix)}(\l^{-1})$$ is a solution for the linear system
\eqref{ay} (Lax pair in $\l$).  Let $(\ti x, \ti y, \ti t)^t=
h(x,y,t)^t$.  A computation gives
$$\ti y+ i\ti x= \frac{i}{\sin\o}
\ (y+  \a\xi +\a^{-1}\eta). $$ Let $\ti \l= h\ast \l$. Then $\ti
\l= \frac{e^c(e^s\l-1)}{e^s\l +1}$, so the pole of this expression
is when $\ti \l=-i$, i.e, when $\frac{e^c(e^s\l-1)}{e^s\l+1}= -i$.
But $r=e^s$ and $e^c= \cot(\o/2)$, so the pole is at  $\l=
\a^{-1}$.  This shows that $h\cdot F$ has one simple pole at $\l=
\a^{-1}$. Note that $h\cdot F$ is equal to
\begin{align*}
(h\cdot F)(x,y,t,\l) &=g_{i,\pi_0(\ti y +i\ti
x)}((\sigma(h)\ast\l)^{-1}) \\
&=\pi_0(\ti y +i\ti x) -e^{i\o} \frac{\l^{-1}-\bar \a}{\l^{-1}-\a}
\pi_0^\perp (\ti y +i\ti x)  \\
&=(\pi(x,y,t)-e^{i\o} \pi^\perp(x,y,t))g_{\a,\pi(x,y,t)}(\tau),
\end{align*}
where $\pi(x,y,t)=\pi_1(y+\a\xi+\a^{-1}\eta)$ and $\pi_1(z)= \pi_0(iz/\sin\o)$ is holomorphic.  
So  $h\cdot F$ is gauge equivalent to the $1$-soliton $g_{\a, \pi}$. Note since $\pi_0:S^2\to \Gr(k, \C^n)$
is smooth and $h\in SO(2,1)$, the monopole given by $g_{\a, \pi}$
decays at spacial variables.
\end{proof}

\bs
\section{B\"acklund transformations and construction of
soliton monopoles\/} \label{da}

Multisolitons with simple poles were constructed by Ward \cite{Wa88}.  Ward, Ioannidou, and Anand (\cite{Wa95, Io96, An98}) derived methods for computing solitons which have poles with higher multiplicities. These multisolitons have dramatic physical properties.  we give here a brief description of a method of ``superposing'' solitons, which is closely related to the permutability formula for B\"acklund transformations.  This technique allowed Dai and Terng \cite{DaTe04} to construct solitons with arbitrary number of poles with arbitrary multiplicities.  

Intuitively, the permutability formula is based on factoring frames.  Given the frames of two solutions $\psi_1$ and $\psi_2$ with singularities at different sets $S_1$ and $S_2$, $S_1\cap S_2=\emptyset$, in $\C\cup\{\infty\}$, we write
$$\psi_3= \ti\psi_1 \psi_2 = \ti\psi_2\psi_1,$$
(i.e., factor $\psi_1\psi_2^{-1}= \ti \psi_2^{-1}\ti \psi_1$).  Here $\psi_3$ has the singularities at $S_1\cup S_2$, and $\psi_j$ and $\ti\psi_j$ have the same singular set $S_j$.  It is not difficult to see that $\psi_3$ is a frame for a solution when $\psi_1$ and $\psi_2$ are.  The details of allowing limiting case where $S_1\to S_2$ yield the interesting but complex solitons.  The converse of factoring solutions is also true, but not completely straightforward.  We now go to the details.

\bthm [{\bf Algebraic B\"{a}cklund transformation\/}] \label{bv}
Let $\psi(x,y,t,\tau)$ be a Ward frame with $P=
(\tau\p_y\psi-\p_\xi\psi)\psi^{-1}$ and $Q=(\tau\p_\eta\psi-\p_y
\psi)\psi^{-1}$, and $g_{\a, \pi}$ a $1$-soliton Ward frame.  Suppose
$\psi$ is holomorphic and non-degenerate at $\tau =\a$. Let
$\ti\pi(x,y,t)$ denote the Hermitian projection of $\C^n$ onto
$\psi(x,y,t,\a)(\Im(\pi(x,y,t))$.  Then 
\ben 
\item $\ti\psi
(x,y,t,\tau)
=g_{\a,\tilde{\pi}(x,y,t)}(\tau)\psi(x,y,t,\tau)g_{\a,\pi(x,y,t)}(\tau)^{-1}$
is holomorphic and non-degenerate at $\tau= \a, \bar \a$, \item
$\psi_1= g_{\a,\ti\pi}\psi = \ti\psi g_{\a,\pi}$ is again a
Ward frame such that
$$
\bca (\tau\p_y\psi_1-\p_\xi \psi_1)\psi_1^{-1}= \ti P,\\
(\tau\p_\eta \psi_1- \p_y \psi_1)\psi_1^{-1}=\ti Q,
\eca
$$
where
$\ti P= P+(\bar{\a}-\a)\p_y\tilde{\pi}$ and $\ti Q= Q +
(\bar{\a}-\a) \p_\eta\tilde{\pi})$.
\een
We will use
$g_{\a,\pi}\ast \psi$ to denote $\psi_1$.
\end{thm}

\begin{proof} We
give a sketch the proof. Statement (1) can be proved by computing
the residue of $\ti \psi$ at $\tau=\a$ and show that it is zero.
We use Proposition \ref{bo} to prove (2).  Set $\ti g= g_{\a,
\ti\pi}$.  Since $D=(\tau\p_y-\p_\xi)$ is a derivation, we
have
$$
(D\psi_1)\psi_1^{-1}= (D\ti g)\ti g^{-1}+ \ti g
(D\psi)\psi^{-1} \ti g^{-1} = (D\ti g)\ti g^{-1}+ \ti g P \ti
g^{-1},
$$
so it is holomorphic for $\tau \in \C\setminus\{\a\}$ and
has a simple pole at $\tau=\a$.  But $\psi_1$ is also equal to
$\ti \psi g$ (here $g=g_{\a,\pi}$ is a $1$-soliton Ward
frame), so
$$(D\psi_1)\psi_1^{-1}= (D\ti \psi)\ti\psi^{-1} + \ti\psi
(Dg)g^{-1}\ti\psi^{-1}.
$$
But $(Dg)g^{-1}$ is independent of $\tau$
and $\ti\psi$ is holomorphic and non-degenerate at $\tau=\a, \bar
\a$, hence the RHS is holomorphic at $\tau=\a$.  So
$(D\psi_1)\psi_1^{-1}$ is holomorphic in $\C$.  But the residue of
$(D\psi_1)\psi_1^{-1}$ at $\tau=\infty$ is also zero.  Hence it
must be independent of $\tau$. Similar argument implies that
$(\tau\p_\eta\psi-\p_x\psi)\psi^{-1}$ is also independent of
$\tau$, so by Proposition \ref{bo}, $\psi_1$ is a Ward frame.

Set $\ti P= (D\psi_1)\psi_1^{-1}= (D\ti g)\ti g^{-1} +\ti g P \ti
g^{-1}$.  Evaluate the residue at $\tau=\infty$ to get $\ti P= \ti
P + (\bar \a-\a)\p_y\ti\pi$.  Similarly, we get the formula for
$\ti Q$.

Since $P, Q$ decay at spacial infinity, $\I-\psi(\cdots, \a)$
decays at spacial infinity.  But $g_{\a, \pi}$ is a $1$-soliton
monopole frame, so $\p_x\pi, \p_\eta\pi$ also decay at spacial
infinity.  Hence $\ti P, \ti Q$ decays at spacial infinity.
\end{proof}

\bs
 \ni{\bf $k$-soliton monopole frames with only simple poles}

\ms
 Let $\a_1, \ldots, \a_k$ be distinct complex numbers and
$\Im(\a_j)>0$ for all $1\leq j\leq k$, $\pi_j^0:S^2\to \Gr(k_j,
\C^n)$ holomorphic maps, and $\pi_j(x,y,t)= \pi_j^0(y+\a_j \xi
+\a_j^{-1}\eta)$. Then $g_{\a_j, \pi_j}$ is a $1$-soliton Ward
frame.  Apply the Algebraic BT (Theorem \ref{bv}) repeatedly as
follows: Set $\psi_1= g_{\a_1, \pi_1}$, and define $\psi_j$
inductively by $\psi_{j}= g_{\a_j, \pi_j}\ast \psi_{j-1}$ for
$2\leq j\leq k$.  Then $\psi_k$ is a $k$-soliton
Ward frame with $k$ simple distinct poles at $\a_1, \ldots,
\a_k$.   These are the same soliton Ward frames constructed by
Ward using the solution to the Riemann-Hilbert problem.

\ms

\ni{\bf $k$-soliton monopole frames with pole data $(\a, k)$}

\ms

Ward's limiting construction is as follows:  Let $f_0, f_1$ be
rational maps from $\C$ to $\C^2$, and $\pi_{1,\e}$ and $\pi_{2,
\e}$ the projections of $\C^2$ onto the complex line spanned by
$f_0(w_{i+\e})+\e f_1(w_{i+\e})$ and $f_0(w_{i-\e})-\e f_1(w_{i-\e})$ respectively, where $w_{i\pm \e}= y+(i\pm\e)\xi +(i\pm\e)^{-1}\eta$. Ward
showed that$$\psi= \lim_{\e\to 0} g_{i-\e, \pi_{2,\e}}\ast
g_{i+\e, \pi_{1,\e}}$$is a $2$-soliton Ward frame
with a double pole at $\tau=i$ and is not stationary.

Since the algebraic BT is easy to compute, Ward's limiting method
can be calculated systematically as follows (for detail see the
paper by Dai and Terng \cite{DaTe04}):  Let $\a_\e=\a+\e$,
$a_j:\C\to \C^n$ be rational maps, and $f_{j,\e}= a_0+ a_1\e +
\ldots + a_{j-1} \e^{j-1}$.  Let $\pi_{j,\e}(x,y,t)$ denote the
Hermitian projection of $\C^n$ onto the complex line spanned by
$f_{j,\e}(y+ \a_\e \xi +\a_\e^{-1}\eta)$.   Set $\psi_1=
\lim_{\e\to 0} g_{\a_\e, \pi_{1,\e}}$, which is a $1$-soliton
frame $g_{\a, \pi}$ (here $\pi$ is the projection onto $\C
a_0(y+\a\xi +\a^{-1}\eta)$). Define $\psi_k$ inductively by
$\psi_k=\lim_{\e\to 0} g_{\a_\e, \pi_{k, \e}}\ast \psi_{k-1}$.
Then $\psi_k$ is a $k$-soliton frame with pole data $(\a, k)$ (i.e., $\psi_k$ has a single pole at $\tau=\a$ with multiplicity $k$).  Note that $\psi_k$ depend on $k$ holomorphic maps from $\C P^1$ to $\cup_{i=1}^{n-1}\Gr(i, \C^n)$.  

\ms

\ni {\bf Soliton frames with arbitrary pole data}

\ms

To get Ward soliton frames with arbitrary pole data, we need  a more
general BT for adding a $k$-soliton with pole data $(\a, k)$ to an
existing Ward  frame (cf. \cite{DaTe04}):

\bthm [{\bf Adding a $k$-soliton with pole data $(\a,
k)$\/}] \label{cp} \par

Suppose $\psi$ is a Ward frame that is
holomorphic and non-degenerate at $\tau=\a, \bar\a$, and $\phi$ a
$k$-soliton monopole frame with pole data $(\a, k)$. Then there
exist unique $\ti\phi$ and $\ti \psi$ such that $\ti\phi \psi=
\ti\psi \phi$, $\ti\phi$ has pole data $(\a, k)$, and $\ti \psi$
is holomorphic and non-degenerate at $\tau= \a, \bar\a$.
Moreover, $\hat \psi= \ti \phi\psi = \ti\psi \phi$ is again a
Ward frame and $\ti \phi$ and $\ti \psi$ are
constructed algebraically.
\ethm

As a consequence,  we see that the
two BTs and the limiting method give rise to Ward soliton frames with arbitrary  pole data.   The following theorem was proved in \cite{DaTe04}. 

\bthm [\cite{DaTe04}] Algebraic BTs, Adding $k$-soliton BTs, and
the limiting method produce all soliton monopoles up to gauge
equivalence.
\ethm

If $\pi_0:S^2\to \Gr(m, \C^n)$ is holomorphic, then the limit of the $1$-soliton frame,
$$\lim_{||(x,y)||\to \infty} \I+ \frac{\a-\bar \a}{\tau-\a}\ \pi_0^\perp(y+\a \xi+\a^{-1}\eta) = \I + \frac{\a-\bar \a}{\tau-\a} \ \pi_0^\perp(\infty)$$
 exists as $(x,y)$ tends to infinity and is independent of $t$. In other words, the Ward soliton frame tends to a fixed rational map $h(\tau)$ at spacial infinity and is independent of time.  It can be checked easily that this property is preserved under the Algebraic BT and limiting  method.  Hence we have
 
 \bprop \label{dh}
 If $\psi$ is a Ward soliton frame, then  
 $$\lim_{||(x,y)||\to \infty} \psi(x,y,t, \tau)$$ exists and is independent of $t$.
 \eprop

\bs
\section{Monopoles with both continuous and discrete scattering
data\/}\label{db}

The Lax pair \eqref{au} of the monopole equation is equivalent to 
\begin{equation}\label{dc}
[\tau(\K_y- \phi)- \K_\xi,\ \tau\K_\eta- \K_y -\phi]=0.
\end{equation}
The linear system associated to this Lax pair is 
\begin{equation}\label{de}
\bca (\tau\p_y -\p_\xi)\psi = (\tau A_y + \tau\phi-A_\xi)\psi,\\
(\tau\p_\eta-\p_y)\psi =(\tau A_\eta+\phi-A_y )\psi.
\eca
\end{equation}

The Algebraic BT theorem for the monopole equation can be proved the same way as for the Ward equation.  We only state the result:

\bthm [{\bf Algebraic BT for Monopoles}]\label{dd}
Suppose $\a\in \C\setminus \R$ is a constant, and $\psi$ is a frame of the monopole solution $(A,\phi)$ (i.e.,  solution of \eqref{de}), and $\psi(x,y,t,\tau)$ is holomorphic and non-degenerate at $\tau= \a$.  Let $g_{\a, \pi}$ be a $1$-soliton Ward frame, $\ti\pi(x,y,t)$ the Hermitian projection onto $\psi(x,y,t,\a)(\Im\pi(x,y,t))$, and
$$\ti\psi= g_{\a, \ti\pi} \psi g_{\a, \pi}^{-1}.$$
Then 
\ben
\item $\ti\psi$ is holomorphic and non-degenerate at $\tau=\a$,
\item $\psi_1= g_{\a, \ti\pi} \psi= \ti\psi g_{\a,\pi}$ is a frame for \eqref{dc} with $\ti A, \ti \phi$ given by 
$$\bca
\ti A_\eta= A_\eta,\\
\ti A_\xi = (1-\frac{\bar a}{\a}) (\p_\xi\ti \pi) h + h^{-1}A_\xi h,\\
\ti A_y +\ti\phi= A_y+\phi,\\
\ti A_y -\ti \phi= (1-\frac{\bar\a}{\a})(\p_y\ti\pi) h + h^{-1}(A_y-\phi) h,
\eca
$$
where $h=\ti\pi +\frac{\a}{\bar\a} \ \ti\pi^\perp$.  
\een
\ethm

Suppose $k:\C\to GL(n,\C)$ is  meromorphic, $k(\infty)=\I$, and  $k(\bar\tau)^*k(\tau)=\I$.  Then
$$\ti k(x,y, t,\tau)= k(y+\tau \xi +\tau^{-1}\eta)$$ satisfies
$$(\tau\p_y-\p_\xi)\ti k=0, \quad (\tau\p_\eta-\p_y)\ti k=0.$$
So if $\psi$ is a solution of \eqref{de} for the monopole solution $(A,\phi)$, then so is  
 $\psi \ti k$.   However, if $\lim_{||(x,y)||\to \infty} \ti k(x,y,t,\tau)$ exists and is independent of $t$, then $k$ must be the constant map $\I$.  Hence we can use this 
 condition to normalize frames to get a unique one:

\bdefn A solution $\psi$ of \eqref{de} for the monopole $(A,\phi)$ is called the {\it normalized monopole frame\/} if 
\ben
\item $\psi(x,y,t,\bar\tau)^*\psi(x,y,t,\tau)=\I$,
\item there exists a map $h(\tau)$ such that $\lim_{||(x,y)||\to \infty} \psi(x,y,t,\tau)= h(\tau)$ exists and is independent of $t$.  
\een
\edefn  

If $\psi$ is a normalized monopole frame, then 
$$E(x,y,t,\mu)= \psi(x,y,t, i(\mu-1)/(\mu+1))$$ is a monopole frame for the Lax pair \eqref{ap} with spectral parameter $\mu$. Moreover,  $E(\cdots, \bar\mu^{-1})^*E(\cdots, \mu)=\I$ and there exists $k(\mu)$ such that $$\lim_{||(x,y)||\to \infty}E(x,y,t,\mu)= k(\mu).$$   We call such $E$ also a normalized monopole frame.  

By Proposition \ref{dh}, a Ward soliton frame is a normalized monopole frame.
 By the Inverse scattering Theorem \ref{bz}, given a smooth map
$s:\R\times S^1\to GL(n,\C)$ such that $\I-s(\cdot, e^{i\o})$
decays for each $\o$ and $s^*=s\geq 0$, then there exists a
solution $E(x,y,t,\mu)$ of the linear system \eqref{ao} such that
$$((E^-)^{-1}E^+)(x,y,0,e^{i\o})= s(x\cos\o +y\sin\o,
e^{i\o})
$$
$E(x,y,t,\mu)\to \I$ as $||(x,y)||\to \infty$,
and  $E$ is holomorphic in $|\mu|\not=1$.   Such $E$ is a normalized monopole frame with only continuous scattering data.

If we apply Algebraic BTs
and General Algebraic BTs  repeatedly to a normalized monopole frame with only continuous scattering data, then we obtain 
normalized monopoles frames with both continuous and discrete
scattering data.  So we get

\bthm  \label{dg}
Let $s:\R\times S^1\to GL(n,\C)$ be a smooth map such that $\I-s(\cdot, e^{i\o})$
decays for each $\o$ and $s^*=s\geq 0$, and $\phi_j$ a normalized soliton monopole frame with pole data $(\a_j, n_j)$ for $j=1, \ldots, k$. Then there is a unique normalized monopole frame $E(x,y,t,\mu)$ such that 
\ben 
\item $E$ is holomorphic for $\mu\in \C\setminus (S^1\cup \{\a_1, \ldots, \a_k\})$, has poles at $\a_j$ with multiplicity $n_j$, and 
$$E_\pm(x,y,t, e^{i\o})=\lim_{|\mu|^{\pm 1}>1, \mu\to e^{i\o}}E(x,y,t,\mu)$$
exist and are smooth,
\item $(E_-^{-1}E_+)(x,y,t,e^{i\o})= s(x\cos\o + y\sin\o -t, e^{i\o})$ is the continuous scattering data of $M$,
\item $E\phi_j^{-1}$ is holomorphic and non-degenerate at $\mu= \a_j$  for $1\leq j\leq k$. 
\een
\ethm

Below we want to prove that all normalized monopole frames with only finitely many poles and a jump across $S^1$ are constructed by the above method.  

First we need to recall a factorization result.  Let $\cd$ denote the group of $f:S^2=\C\cup \{\infty\}\to GL(n,\C)$ satisfies the following conditions:
\ben
\item $f(\bar\mu^{-1})^*f(\mu)=\I$,
\item $f_\pm(e^{i\o}):=\lim_{|\mu|^{\pm 1}>1, \mu\to e^{i\o}} f(\mu)$ exist and are smooth,
\item  $f$ is holomorphic in $S^2\setminus S^1$ except  with possible finitely many poles.
\een

The following results were proved in \cite{TeUh00}:
\ben
\item[(a)] Given $f\in \cd$, then there exist uniquely $h_i, g_i\in \cd$ such that $f=h_1g_1= g_2 h_2$, $f_i$ is holomorphic in $S^2\setminus S^1$, $g_i$ are rational maps, and $g_i(-1)=\I$  for $i=1, 2$. 
\item[(b)] If $f\in \cd$  has a simple pole at $\a$, then there is a
 unique projection $\pi$ such that $fg_{\a, \pi}^{-1}$ is holomorphic and non-degenerate at $\a$.
 \een

Note that if $E$ is a normalized monopole frame then $E(x,y,t,\cdot)\in \cd$. Moreover, 
\ben 
\item[(i)] $E$ is a soliton frame if and only if $E$ is rational in $\mu$, 
\item[(ii)] $E$ has only continuous  scattering data if and only if  $E(x,y,t,\cdot)$ is holomorphic in $S^2\setminus S^1$.  
\een
For general normalized monopole frames with both poles and jumps across $S^1$, we have

\bthm  [{\bf Subtracting $1$-soliton\/}]  \label{by}
Suppose $\psi$ is
a normalized monopole frame, and $\psi$ has a simple pole at
$\tau=\a$ (may have other singularities as well).  Then there
exist unique $\ti\psi, \psi_1$ and smooth $\pi, \ti\pi:\R^{2,1}\to
\Gr(k, \C^n)$ such that
 \ben
 \item[(i)] $\psi= \ti\psi g_{\a, \pi} =
g_{\a, \ti\pi} \psi_1$,
 \item[(ii)] $\ti\psi$ and $\psi_1$ are holomorphic and
non-degenerate at $\tau=\a, \bar \a$,
\item[(iii)] $g_{\a,\pi}$ is a
normalized $1$-soliton monopole frame, and $\psi_1$ is a
normalized monopole frame. \een \ethm

\begin{proof} Statements (i) and (ii) follows from known results stated just before the theorem.  To prove (iii) we use residue calculus.  Set $D=\tau\p_y-\p_\xi$, and let
$g:=g_{\a,\pi}$, $\ti g:=g_{\a,\ti\pi}$. Then
$$(D\psi)\psi^{-1}= (D\ti \psi)\ti\psi^{-1}+\ti\psi
(Dg)g^{-1}\ti\psi.
$$
Since the LHS is a degree one polynomial in $\tau$, it is holomorphic at $\tau=\a$.  So the
residue of the RHS at $\tau=\a$ must be zero, which implies that
$$
-\ti\psi(\cdots,
\a)(\a\p_y\pi -\p_\xi \pi)\pi \ti\psi(\cdots, \a)^{-1}=0.
$$
But $\ti\psi(\cdots, \a)$ is non-degenerate, hence
$(\a\p_y\pi-\p_\xi \pi)\pi=0$.  Similarly, calculate the residue
at $\tau=\a$ in the expression $(\tau\p_\eta\psi-
\p_y\psi)\psi^{-1}$ to get $(\a\p_\eta\pi-\p_y\pi)\pi=0$.  By
Proposition \ref{bw}, $g_{\a, \pi}$ is a $1$-soliton monopole
frame. Since $\psi_1= \ti g^{-1}\psi$,
$$(D\psi_1)\psi_1^{-1}= -\ti
g^{-1}D\ti g + \ti g^{-1}(D\psi)\psi^{-1} \ti g.
$$
The LHS is holomorphic at $\tau=\a, \bar\a$ and the RHS is
holomorphic for all $\tau\not=\a, \bar\a$. So
$(D\psi_1)\psi_1^{-1}$ is a degree one polynomial in $\tau$.  Similarly,
$(\tau\p_\eta \psi_1-\p_x\psi_1)\psi_1^{-1}$ is a degree one polynomial  in
$\tau$. By Proposition \ref{bo}, $\psi_1$ is a monopole frame.
\end{proof}

Similar argument gives

\bthm\label{cb}  {\bf [Subtracting a soliton with pole data $(\a,
k)$]}

\par
Suppose $\psi$ is a normalized monopole frame, and $\psi$ has a
pole at $\tau=\a$ with multiplicity $k$ (may have other
singularities as well). Then there exist unique $k$-soliton
monopole frame $g$ with pole data $(\a, k)$, a normalized monopole
frame $\psi_1$, and maps $\ti \psi$ and $\ti g$ such that $\psi_1$
and $\ti\psi$ are holomorphic and non-degenerate at $\tau=\a,
\bar\a$ and $\psi= \ti\psi g = \ti g \psi_1$. \ethm

 A consequence of the above two theorems is that every normalized monopole frame with continuous scattering data and finitely many poles can be obtained from Theorem \ref{dg}.

\bs

\section{Appendix: Existence of Continuous Scattering Data for Small Solutions\/}

The over-all details of the scattering data described in sections 7
 and 10
where we find solutions with a combination of a jump across the
unit circle and point (singularities) measures does fit into the
general scheme proposed by Beals and Coifman  \cite{BeaCoi85b}, \cite{BeaCoi89}. 
Beals and Coifman point out that the Self-dual Yang-Mills equations
have local scattering data more like the AKNS scattering problem than
either type of KP, and Fokas and Ioannidou point out that the Ward equation
inherits this similarity. Existence of solutions which have combinations
of the two types of scattering data follows from the existence of the
inverse scattering transform and the process of adding discrete measures
by the described B\"acklund transforms.  Papers  of Fokas and Ioannidou \cite{FoIo01} and
Villarroel \cite{Vi90} discuss the existence of the scattering and inverse scattering
transforms. 

We find that, as in the AKNS models, small data leads to scattering data
which consists only of the continuous part.  For simplicity, we assume
all the data  is rapidly decaying and lies in the Schwartz space. Recall that 
The spacial operator is 
$$D_s(\mu) = D_1(\mu)-D_2(\mu) =
\frac{(\mu-\mu^{-1})}{2}\K_x -i\frac{(\mu+\mu^{-1})}{2}\K_y -
i\phi.$$

\bthm
 Assume that there exists a gauge transformation such that in the
given gauge $(A,\phi)$ is small in $W^{2,1}$. Then the chern vector $c(\mu)$ is
zero for all $\mu \in C \cup \{\infty\} \setminus S^1$.  Furthermore there exists continuous
scattering data $S_\o = I + G_\o$ where $G_\o$ is small in
 $L^\infty_{(x,y,t)}$.
\ethm

\bcor
  Given initial data $(A,\phi)$ which is small in $W^{2,1}$, there
exist global in time solutions of the space-time monopole problem with
this initial data.  These solutions are unique up to gauge transformation.
  \ecor

First we explain  how the corollary follows from the theorem.
The existence of global solutions  follows from the existence of
scattering data which is purely continuous for
the initial data, the known  flow of the scattering data
under time, and the existence of the inverse scattering transform. Because
the scattering data remains small, the solution remains small in time.
Suppose there is a time $t_\o$ at which uniqueness fails.  Since we presume
the solution to be continuous in time, the second solution is small for
a short time near $t_\o$.  Hence it has scattering data, and must be identified
with the solution constructed by inverse scattering data by a gauge
transformation.

To prove the theorem, we note that the proof is rather standard away from
the unit circle, although it will follow from the proof we give near the
unit circle as well.  The difficulty is to prove that the limits exist
as $|\mu| \to 1$.  To do this, use the combination of Lorentz and
fractional linear transformation so that $\mu$ is pure imaginary. Now let
 $$\tau = (\mu + {\mu}^{-1})/(\mu - {\mu}^{-1}),$$
so $\tau < 1$ is real.  We assume $|\mu|<1$, so $\tau<0$.  For $|\mu|>1$, then the reality condition implies that  $E_\mu = ((E_{\bar\mu^{-1}})^*)^{-1}$.   The equation for the frame now reads:
 $$(\p_x - i \tau \p_y) E_{\tau} =
 G_{\tau}$$
where $$G_\tau(x,y,t) = A_\tau(x,y,t) E_\tau (x,y,t), \quad A_{\tau} =  A_x - i \tau A_y(x,t) -\sqrt{1 - \tau^2}\ \Phi.$$

As is usual with solutions with small data, we solve by iteration. 
Let  $E_\tau = \I + \sum_ {j\geq 1} Q_{j,\tau}$. Set $Q_{0, \tau}=\I$,  $Q_{j,\tau} (-\infty, y) = 0$, and
define $Q_{j,\tau}$ iteratively by
    $$(\p_x - i \tau \p_y)Q_{j,\tau}
 =  G_{j-1,\tau} = A_\tau Q_{j-1,\tau}.$$
Here $A_{\tau}$ is as above.

   We now take the Fourier transform in the $y$ variable alone, and denote 
the transform in $y$ only of an expression $Q$  by $\hat Q$.  Let $\zeta$ be the transform variable of $y$.  We have
$$ \p_x \hat Q_{j,\tau} - \zeta{\tau}\hat Q_{j,\tau}
=\hat G_{j-1,\tau}.$$ 
  But this ODE can be solved explicitly:
  $$\hat Q_{j,\tau}(x,\zeta, t)=\sgn(\zeta\tau) \int_{\zeta\tau(x-r)\leq 0} e^{\zeta\tau(x-r))} \hat G_{j-1, \tau}(r, \zeta, t)\ dr.$$
  So 
  $$\int_{-\infty}^\infty \hat Q_{j, \tau}(x,\zeta,t) d\zeta = \int_{-\infty}^\infty \sgn(\zeta\tau)
  \int_{\zeta\tau(x-r)\leq 0} e^{\zeta\tau(x-r))} \hat G_{j-1, \tau}(r, \zeta, t)\ dr \ d\zeta.$$
   If $\tau \zeta>0$, then the above integral is equal to 
  $$\int_{-\infty}^\infty\int_{-\infty}^0 e^{\tau \zeta s} \hat G_{j-1, \tau}(x-s, \zeta, t)\ ds d\zeta.$$ 
 For $\tau \zeta<0$, we get a similar but different formula. This explains the jump. 
  So we  obtain in a straightforward fashion the estimate independent of $\tau$:
 $$|| \hat Q_{j,\tau}||_1 \leq  || \hat G_{j-1,\tau}||_2 =
  ||\hat A_{\tau}\ast \hat Q_{j-1,\tau}||_2 \leq
       ||\hat A_{\tau}||_2 ||\hat Q_{j-1,\tau}||_1.$$
The norms used are $||\cdot  ||_1$ as the  $L^\infty_x L^1_{\zeta}$ and $||\cdot  ||_2$
as the $L^1_{(x,\zeta)}$ norm.  Given these estimates, the solution 
$$\hat Q(x,\zeta,t) =  \sum_{j\geq 0} \hat Q_j(x,\zeta,t)$$
 exists in $L^{\infty}_xL^1_{\zeta}$  when  $||\hat A_{\tau}||_2 \leq 1$.  But it is
easy to see that when the $W^{2,1}$ norm of $(A, \phi)$ is small, this inequality is
satisfied.  Furthermore, $||\hat Q||_1$ bounds $||Q||_{\infty}$ and existence
follows.  Iterating these same estimates in the various derivatives of
$Q$ will give estimates on the higher norms, although we do not expect or
require smallness in the higher derivatives and weights.

\bs


\begin{thebibliography}{99}

\bibitem{An97}
Anand, C. K., \emph{{W}ard's solitons},  Geom. Topol., {\textbf 1}
(1997), 9--20.
%
\bibitem{An98}
Anand, C. K., \emph{{W}ard's solitons II, Exact solutions}, Canad.
J. Math., {\textbf 50} (1998), 1119--1137.
%
\bibitem{AtHi88}
Atiyah, M. F. and Hitchin, N. J., \emph{{T}he geometry and
dynamics of magnetic monopoles}, Princeton University Press,
Princeton, New Jersey, 1988.
%
\bibitem{BeaCoi85b} 
Beals, R., Coifman, R.R., \emph{{M}ultidimensional inverse scattering and nonlinear partial differential equations}, Proc. Symp. Pure Math., \textbf{43} (1985), 45-70
%
\bibitem{BeaCoi89}
Beals, R., Coifman, R.R., \emph{{L}inear spectral problems, non-linear
equations and the $\bar\partial$-method}, Inverse Problems, \textbf{5} (1989), 87-130
%
\bibitem{DaTe04}
Dai, B. and Terng, C. L., \emph{{B}\"acklund transformation, Ward
solitons, and unitons}, arXiv:math.DG/0405363.
%
\bibitem{FoIo01}
Fokas, A. S. and Ioannidou, T. A., \emph{{T}he inverse spectral
theory for the Ward equation and for the $2+1$ chiral model},
Comm. Appl. Anal., {\textbf 5} (2001), no. 2, 235--246.
%
\bibitem{Hi83}
Hitchin, N. J., \emph{{O}n the construction of monopoles}, Comm.
Math. Phys., {\textbf 89} (1983), 145--190.
%
\bibitem{Hi87}
Hitchin, N. J., \emph{{T}he self-duality equations on a Riemann
surface}, Proc. London Math. Soc., {\textbf 55} (1987), 59--126.
%
\bibitem{Io96}
Ioannidou, T., \emph{{S}oliton solutions and nontrivial scattering
in an integrable chiral model in $(2+1)$ dimensions}, J. Math.
Phys., {\textbf 37} (1996), 3422--3441.
%
\bibitem{IoZa98a}
Ioannidou, T. and Zakrzewski, W., \emph{{S}olutions of the
modified chiral model in $(2+1)$ dimensions}, J. Math. Phys.,
{\textbf 39} (1998) no.5, 2693--2701.
%
\bibitem{IoZa98b}
Ioannidou, T. and Zakrzewski, W., \emph{{L}agrangian formulation
of the general modified chiral model}, Phys. Lett. A, {\textbf
249} (1998), 303--306.
%
\bibitem{MaZa81}
Manakov, S. V. and Zakharov, V. E., \emph{{T}hree-dimensional
model of relativistic-invariant theory, integrable by the inverse
scattering transform}, Lett. Math. Phys., {\textbf 5} (1981),
247--253.
%
\bibitem{PrSe86}
Pressley, A. and Segal, G., \emph{{L}oop groups}, Oxford
University Press, 1986.
%
\bibitem{Tao01}
Tao, T.,  \emph{{G}lobal regularity of wave maps. II.  Small
energy in two dimensions}, Comm. Math. Phys, {\textbf 224} (2001),
no. 2, 443--544.
%
\bibitem{Tat04}
Tataru, D.,  \emph{{T}he wave maps equation}, Bull. Amer. Math.
Soc., {\textbf 41} (2004), no. 2, 185--204.
%
\bibitem{TeUh98}
Terng, C. L. and Uhlenbeck, K.,  \emph{{P}oisson actions and
scattering theory for integrable systems}. Surveys in differential
geometry: integrable systems, 315--402, Surv. Diff. Geom., IV,
International Press, Boston, MA, 1998.
%
\bibitem{TeUh00}
Terng, C. L. and Uhlenbeck, K.,  \emph{{B}\"{a}cklund
transformations and loop group actions}, Comm. Pure Appl. Math.,
{\textbf 53} (2000), 1--75.
%
\bibitem{Uh89}
Uhlenbeck, K., \emph{{H}armonic maps into Lie groups (classical
solutions of the chiral model)}, J. Diff. Geom., {\textbf 30}
(1989), 1--50.
%
\bibitem{Uhl92}
Uhlenbeck, K., \emph{{O}n the connection between harmonic maps and the self-dual Yang-Mills and the sine-Gordon equations}, J. Geom. Phys. \textbf{8} (1992), 283-316
%
\bibitem{Vi90}
Villarroel, J., \emph{{T}he inverse problem for Ward's system},
Stud. Appl. Math., {\textbf 83} (1990), 211--222.
%
\bibitem{Wa88}
Ward, R.S., \emph{{S}oliton solutions in an integrable chiral
model in $2+1$ dimensions}, J. Math. Phys., {\textbf 29} (1988),
386--389.
%
\bibitem{Wa90}
Ward, R.S., \emph{{C}lassical solutions of the chiral model,
unitons, and holomorphic vector bundles}, Comm. Math. Phys.,
{\textbf 128} (1990), 319--332.
%
\bibitem{Wa95}
Ward, R.S., \emph{{N}ontrivial scattering of localized solutions
in a $(2+1)$-dimensional integrable systems}, Phys. Letter A,
{\textbf 208} (1995), 203--208.
%
\bibitem{Wa99}
Ward, R. S. \emph{{I}ntegrable systems and twistors}, in Integrable Systems, Oxf. Grad. Texts Math., {\textbf 4}  (1999), 121-134.
%
\bibitem{ZaMi78}
Zakharov, V. E. and Mikhailov, A. V., \emph{{R}elativistically
invariant two dimensional models of fields theory which are
integrable by means of the inverse scattering problem method},
Sov. Phys. JETP, {\textbf 47} (1978), no. 6, 1017--1027.
%

\end{thebibliography}
\end{document}